\documentclass[12pt,reqno]{amsart}
\usepackage{graphicx, amssymb, xcolor,pdfsync}
\usepackage[all,cmtip]{xy}
\numberwithin{equation}{section}
\usepackage{mathrsfs}
 \usepackage{mathpazo}
\usepackage{hyperref}
\usepackage{graphicx}
\usepackage{yhmath}
\usepackage{mathdots}
\usepackage{MnSymbol}

\hypersetup{
    colorlinks=true,       
    linkcolor=blue,          
    citecolor=blue,        
    filecolor=blue,      
    urlcolor=blue           
}
\usepackage{tikz}
\usetikzlibrary{matrix,arrows}
\DeclareMathOperator{\Aut}{Aut}
\usepackage[switch]{lineno}
\usepackage{yfonts}

\advance\hoffset by -0.9 truecm
\DeclareMathOperator{\Pm}{Prym}
\DeclareMathOperator{\exppp}{exp}

\begin{document}
\newcommand{\s}{\vspace{0.2cm}}
\newcommand{\Sp}{\mbox{Sp}}
\newcommand{\sy}{\mbox{sym}}

\newtheorem{theo}{Theorem}
\newtheorem{prop}{Proposition}
\newtheorem{coro}{Corollary}
\newtheorem{lemm}{Lemma}
\newtheorem{claim}{Claim}
\newtheorem{example}{Example}
\theoremstyle{remark}
\newtheorem*{rema}{\it Remark}
\newtheorem*{rema1}{\it Remark}
\newtheorem*{defi}{\it Definition}
\newtheorem*{theo*}{\bf Theorem}
\newtheorem*{coro*}{Corollary}
\newtheorem{hecho}[theo]{\bf Hecho}

\title[Abelian varieties and Riemann surfaces with group action]{Abelian varieties and Riemann surfaces with \\generalized quaternion group action}
\date{}

\author{Angel Carocca}
\address{Departamento de Matem\'atica y Estad\'istica, Universidad de La Frontera, Avenida Francisco Salazar 01145, Temuco, Chile.}
\email{angel.carocca@ufrontera.cl}

\author{Sebasti\'an Reyes-Carocca}
\address{Departamento de Matem\'atica y Estad\'istica, Universidad de La Frontera, Avenida Francisco Salazar 01145, Temuco, Chile.}
\email{sebastian.reyes@ufrontera.cl}

\author{Rub\'i E. Rodr\'iguez}
\address{Departamento de Matem\'atica y Estad\'istica, Universidad de La Frontera, Avenida Francisco Salazar 01145, Temuco, Chile.}
\email{rubi.rodriguez@ufrontera.cl}

\thanks{Partially supported by Fondecyt Grants 1200608, 11180024 and 1190991}
\keywords{Riemann surfaces, Abelian varieties, group actions, automorphisms}
\subjclass[2010]{30F10, 14H37, 14H40, 30F20, 14K02, 14H30}

\begin{abstract} In this article we consider  Riemann surfaces and  abelian varieties endowed with a group of automorphisms isomorphic to a generalized quaternion group. We provide isogeny decompositions of each abelian variety with this action, compute dimensions of the corresponding factors and provide conditions under which this decomposition is nontrivial. We then specialize our results to the case of Jacobians and relate them to the so-called genus-zero actions on  Riemann surfaces. We also give a complete classification and description of the complex one-dimensional families of Riemann surfaces and Jacobians with a generalized quaternion group action, extending known results concerning the quasiplatonic case. Finally, we construct  and describe explicit families of abelian varieties with a quaternion group action and derive a period matrix for the Jacobian of the surface with  full automorphism group of second largest order among the hyperelliptic surfaces of genus four.

\end{abstract}
\maketitle
\thispagestyle{empty}

\section{Introduction and statement of the results}

Complex abelian varieties and compact Riemann surfaces have been extensively studied since the nineteenth century and their foundations go back to works of Abel, Jacobi, Riemann, Klein and Hurwitz among others. These objects, that are different generalizations of the classically known elliptic curves, are intimately related to each other by the Torelli map which takes a surface and sends it to its  Jacobian variety. 

\s

Although the theory of groups actions on compact Riemann surfaces and complex abelian varieties is classical and have been widely investigated over the last decades, it still attracts much attention because of several interesting questions in this subject remain unsolved up to the present. 

\s

A foundational result in the area claims that each abstract finite group can be realized as a group of automorphisms of some compact Riemann surface, hence of some Jacobian variety, and therefore of some  complex abelian variety. Articles aimed at studying  group actions of special classes of groups can be found in the literature in plentiful supply.

\s

This article is devoted to considering complex abelian varieties and compact Riemann surfaces possessing a group of automorphisms isomorphic to the generalized quaternion group. These groups form one of the four series of non-abelian $2$-groups that have a cyclic subgroup of index two; the simplest member of the series being the  classical quaternion group of order eight. Throughout the paper, we  denote by $$Q(2^n)=\langle x,y : x^{2^{n-1}}, y^2x^{2^{n-2}}, yxy^{-1}x \rangle $$ the generalized quaternion group of order $2^n,$ where $n \geqslant 3$ is an integer.

\subsection{Abelian varieties with $Q(2^n)$-action} 
Every action of a finite group $G$ on an abelian variety $A$ induces $\mathbb{Q}$-algebra homomorphism $$ \Xi : \mathbb{Q}[G] \to \mbox{End}_{\mathbb{Q}}(A)=\mbox{End}(A) \otimes_{\mathbb{Z}} \mathbb{Q}$$from the rational group algebra of $G$ to the rational endomorphism algebra of $A.$ As proved by Lange and Recillas in \cite{LR}, the decomposition of $1$ as a sum of appropriate idempotents of $\mathbb{Q}[G]$ induces, by means of $\Xi$, an isogeny decomposition $$A \sim A_1 \times \cdots \times A_r$$where the factors $A_i$ are abelian subvarieties of $A$ that are pairwise non-$G$-isogenous. This decomposition is called the {\it isotypical decomposition} of $A$ with respect to $G$. Moreover, each $A_i$ decomposes further isogenously as a suitable power of an abelian subvariety $B_i$ of it, and therefore the following isogeny is obtained $$A \sim B_1^{n_1} \times \cdots \times B_r^{n_r}.$$ This decomposition is called a {\it group algebra decomposition} of $A$ with respect to $G.$ The  subvarieties $A_i$ and $B_i$ together with the integers $n_i$  are closely related to the representations of $G$, as we shall explain later. 
\s

The first result of this paper provides the above-mentioned isogeny decompositions of each abelian variety endowed with a group of automorphisms isomorphic to ${Q}(2^n).$  In order to state it, we need to introduce some notations. 

\s

Let $A$ be an abelian variety with $G$-action. For each subgroup $K$ of $G$ we write $$p_K=\Sigma_{k \in K}k \in \mathbb{Q}[G] \, \mbox{ and } \, A_K = \mbox{Im}(\Xi(p_K)).$$

Observe that  $A_K$ is the abelian subvariety of $A$ on which $K$ acts trivially. In addition, if $K_1 \leqslant K_2$ are subgroups of $G$ then $A_{K_2} \subseteq A_{K_1}$ and, as proved in \cite{CR-H}, there exists an abelian subvariety $P$ of $A_{K_1}$ in such a way that $$A_{K_1} \sim P \times A_{K_2}.$$The abelian subvariety $P$ as before will be denoted by $\Pm(A_{K_1} / A_{K_2}).$ 
\s

Along the article, we shall employ the following notation for the proper nontrivial subgroups of $Q(2^{n}):$ 
$${H}_j=\langle x^{2^{n-j}}, y \rangle \hspace{0.4 cm} K_i=\langle x^{2^{n-i}} \rangle\hspace{0.4 cm} \tilde{H}_j=\langle x^{2^{n-j}}, xy \rangle$$where $i \in \{2, \ldots,n\}$ and $j \in \{2, \ldots, n-1\}.$ Among these subgroups, there are four that will play a special role in our results; we shall use a different notation for them: $$Z:=K_2=\langle y^2 \rangle \hspace{0.4 cm} N_1:=K_n=\langle x \rangle \hspace{0.4 cm} N_2:=H_{n-1}=\langle x^2, y\rangle \hspace{0.4 cm} N_3:=\tilde{H}_{n-1}=\langle x^2, xy\rangle.$$

\begin{theo} \label{isog-abel-gen} Let $n \geqslant 4$ be an integer and 
let $A$ be an abelian variety endowed with a group of automorphisms $G$ isomorphic to $Q(2^n).$ 

The isotypical decomposition of $A$ with respect to $G$ is given by
$$
A \sim A_{G} \times \Pi_{i=1}^3   \Pm(A_{N_i} /A_{G}) \times \Pm(A/ A_{Z}) \times \Pi_{i=2}^{n-2} \Pm(A_{K_i} / A_{K_{i+1}}).
$$ Moreover, two group algebra decompositions of $A$ with respect to $G$ are given by\begin{equation*}A \sim A_{G} \times \Pi_{i=1}^3   \Pm(A_{N_i} /A_{G}) \times \Pm(A/ A_{Z}) \times \Pi_{j=2}^{n-2} \Pm(A_{H_j} / A_{H_{j+1}})^2\end{equation*}\begin{equation*}\label{os2}A \sim A_{G} \times \Pi_{i=1}^3   \Pm(A_{N_i} /A_{G}) \times \Pm(A/ A_{Z}) \times \Pi_{j=2}^{n-2} \Pm(A_{\tilde{H}_j} / A_{\tilde{H}_{j+1}})^2.\end{equation*}
\end{theo}

The previous theorem tells us that each factor arising in the isotypical decomposition and in a group algebra decomposition of $A$ with respect to $G$ is isogenous to an abelian subvariety of $A$ with a geometric meaning. It is worth pointing out that, in general, the aforementioned fact is not true.  Observe, moreover, that a direct consequence of the previous result is the fact that there are isogenies $$\Pi_{i=2}^{n-2} \Pm(A_{K_i} / A_{K_{i+1}}) \sim \Pi_{j=2}^{n-2} \Pm(A_{\tilde{H}_j} / A_{\tilde{H}_{j+1}})^2 \sim \Pi_{j=2}^{n-2} \Pm(A_{{H}_j} / A_{{H}_{j+1}})^2.$$

It will be a consequence of the proof of the theorem that, indeed, there are isogenies $$ \Pm(A_{K_j} / A_{K_{j+1}}) \sim  \Pm(A_{\tilde{H}_j} / A_{\tilde{H}_{j+1}})^2 \sim  \Pm(A_{{H}_j} / A_{{H}_{j+1}})^2$$for each $j \in \{2, \ldots, n-2\}.$

\s

The group algebra decompositions only depend on the algebraic structure of the group. By contrast, the dimensions of the factors do depend on the way the group acts, and this dependence is encoded in the analytic representation of $G$, denoted by $$\rho_a : G \to \mbox{GL}(V) \, \mbox{ where } \, A=V/\Lambda.$$

To state the following proposition we need to fix some notations concerning the well-known  representations of the generalized quaternion group. Hereafter, the complex irreducible representations of  $Q(2^n)$ of degree one will be denote by 
$$\chi_1 : x, y \mapsto 1, \,\, \chi_2: x \mapsto 1, y \mapsto -1, \,\, \chi_3: x \mapsto -1, y \mapsto 1 \,\, \mbox{ and }\,\,  \chi_4: x,y \mapsto -1,$$whereas the ones of degree two will be denoted by  \begin{equation*}
\Theta_s  : x \mapsto 
 \mbox{diag} (\omega^s, \bar{\omega}^s) \, \mbox{ and }\, y \mapsto \left( \begin{smallmatrix}
0 & (-1)^s \\
1 & 0 \\
\end{smallmatrix} \right)
\end{equation*}where $\omega:=\mbox{exp}(2 \pi i / 2^{n-1})$ and $s \in  \{1, \ldots, 2^{n-2}-1\}$.

\begin{prop} \label{pp1} Let $n \geqslant 4$ be an integer and let $A$ be an abelian variety with a group of automorphisms $G$ isomorphic to  $Q(2^n).$ Assume that the analytic representation of the action of $G$ on $A$ decomposes as a sum of irreducible representations as \begin{equation*}\label{imp}\rho_a \simeq \oplus_{j=1}^4 a_j \chi_j \oplus \oplus_{s=1}^{2^{n-2}-1}b_s\Theta_s\end{equation*}for some  non-negative integers $a_1, \ldots, a_4, b_1, \ldots, b_{2^{n-2}-1}.$  

For each $j \in \{2, \ldots, n-1\}$, choose $t_j \in \mathbb{Z}_{2^{n-1}}$  such  that $\omega \mapsto  \omega^{t_j}$  generates the Galois group associated to the character field of $\Theta_{2^{j-1}}$ (which is isomorphic to the cyclic group of order ${2^{n-j-2}}$). Then $$b_1=b_3=\cdots =b_{2^{n-2}-1} \,\, \mbox{ and }\,\, b_{2^{j-1}}=b_{t_j^k2^{j-1}} \, \mbox{ for each } 0 \leqslant k \leqslant 2^{n-j-2}-1,$$ and the dimension of the factors arising in the group algebra decompositions of $A$ with respect to $G$ of Theorem \ref{isog-abel-gen} is given in the following table.

\s

\begin{center}
\begin{tabular}{|c|c|c|c|c|}  
\hline
$A_{G}$  &  $\Pm(A_{N_i}/A_{G})$ & $\Pm(A/A_{Z})$ & $\Pm(A_{H_j} / A_{H_{j+1}})$ & $\Pm(A_{\tilde{H}_j} / A_{\tilde{H}_{j+1}})$\\ \hline 
$a_1$ & $a_{i+1}$ & $2^{n-2}b_1$ & $2^{n-j-2}b_{2^{j-1}}$& $2^{n-j-2}b_{2^{j-1}}$\\ \hline
\end{tabular}
\end{center}

\end{prop}

\s

We emphasize the fact that the dimension of any factor arising in a group algebra decomposition might be equal to zero. Indeed, in an extreme situation all but one of the factors can be zero-dimensional; in this case we say that such a group algebra decomposition is {\it trivial}. Clearly, if a group algebra decomposition is trivial then each group algebra decomposition  is also as well.

\begin{theo} \label{equiv} Let $n \geqslant 4$ be an integer and let $A$ be an abelian variety with a group of automorphisms $G$ isomorphic to  $Q(2^n).$ Assume the analytic representation $\rho_a$ of $G$ to be as in Proposition \ref{pp1}. The following statements are equivalent.
\begin{enumerate}
\item A group algebra decomposition of $A$ with respect to $G$ is trivial.
\item The dimension of $A_Z$ is zero.
\item The dimension of $A_{K}$ is zero, for each nontrivial subgroup $K$ of $G.$
\item The integers $a_i$ and $b_s$ equal zero, for all $i$ and for all even $s.$
\item 1 is not an eigenvalue of $\rho_a(g)$ for each nontrivial $g \in G.$
\end{enumerate}
\end{theo}

A complex representation of a finite group satisfying the last statement of the preceding theorem is called {\it fixed point free}. It would be an interesting problem to classify finite groups acting  on abelian varieties in such a way that its analytic representation is fixed point free; this general problem will be considered elsewhere.

\subsection{Jacobian varieties with $Q(2^n)$-action} Let $S$ be a compact Riemann surface of genus $g \geqslant 2$. The Jacobian variety $JS$ of $S$ is an irreducible principally polarized abelian variety of dimension $g.$ We recall  that, by the classical Torelli's theorem, $$S \cong S' \,\, \iff \,\, JS \cong JS'.$$

Let $G$ be a finite group acting on $S$ and let $K$ be a subgroup of $G.$ The  covering map $\pi_K:S \to S_K$ given by the action of $K$ on $S$ induces a homomorphism of abelian varieties $$\pi_K^*: JS_{K} \to JS \,\, \mbox{ and therefore }\,\, JS \supset  \pi_K^*(JS_{K})\sim JS_{K}.$$Observe that if $A=JS$ then $$A_K=\pi_K^*(JS_K) \sim JS_K \,\, \mbox{ and } \Pm(A_{K_1}/A_{K_2})=\Pm(S_{K_1} \to S_{K_2}),$$where the latter symbol stands for the generalized Prym variety associated to a covering map of compact Riemann surfaces. 

It then follows immediately  that Theorems \ref{isog-abel-gen} and  \ref{equiv} can be restated as follows.

\begin{coro} \label{isog-sup-gen} Let $n \geqslant 4$ be an integer and let $S$ be a compact Riemann surface endowed with a group of automorphisms $G$ isomorphic to $Q(2^n).$ Then a group algebra decomposition of $JS$ with respect to $G$ is given by$$JS \sim JS_{G} \times \Pi_{i=1}^3   \Pm(S_{N_i} \to S_{G}) \times  \Pm(S \to S_{Z}) \times \Pi_{j=2}^{n-2} \Pm(S_{H_j} \to S_{H_{j+1}})^2.$$

If, in addition, the quotient $S_{G}$ is isomorphic to the projective line then $$JS \sim JS_{N_1}  \times  JS_{N_2}  \times JS_{N_3}  \times \Pm(S \to  S_{Z}) \times \Pi_{j=2}^{n-2} \Pm(S_{H_j} \to S_{H_{j+1}})^2.$$

Moreover, the following statements are equivalent.
\begin{enumerate}
\item A group algebra decomposition of $JS$ with respect to $G$ is trivial.
\item The genus of $S_Z$ is zero.
\item The genus of $S_K$ is zero, for each nontrivial subgroup $K$ of $G.$
\end{enumerate}
\end{coro}

Observe that in the case of Jacobians with $Q(2^n)$-action, the dimension of the factors arising in a group algebra decomposition of $JS$ with respect to $G$ can be computed  in terms of the genus of the quotient Riemann  surfaces given by the action on $S$ of certain subgroups of $G.$ See also \cite[Theorem 5.12]{anita-ibero}.


\subsection{Genus-zero actions on Riemann surfaces} Following \cite[Definition 1]{KS}, the action of a  group $G$ on a compact Riemann surface $S$ is termed a {\it genus-zero} action if
$$1_G \neq K \leqslant G \,\, \implies \,\, S_K \cong \mathbb{P}^1.$$

Kallel and Sjerve in \cite{KS} provided a classification of genus-zero actions on Riemann surfaces. They succeeded in proving that only a few  groups  arise with this property: the cyclic groups, the generalized quaternion groups, the polyhedral groups and the so-called ZM groups $C_p \rtimes_2 C_4$ where $p$ is an odd prime.  Furthermore, they gave a topological description by providing the so-called signatures that realize each such actions. 

\s

Let $K$ be a group of automorphisms of $S.$ We recall that the tuple $$\sigma=(\gamma; k_1, \ldots, k_l) \in \mathbb{Z}^{l+1}  \mbox{ where }  \gamma \geqslant 0 \mbox{ and } k_i \geqslant 2$$is called the {\it signature} of the action of  $K$ on  $S$ if the genus of $S_K$ is $\gamma$ and the associated branched regular covering map $S \to S_K$  ramifies over exactly $l$ values $y_1, \ldots, y_l$ and the fibre over $y_i$ consists of points with $K$-isotropy group of order $k_i$ for each $i.$ If $\gamma=0$ and $l=3$ then  we say that such an action is {\it triangle}.

\s

In \cite[Theorem 3]{KS}, it is claimed that the genus-zero actions of $Q(2^n)$ have signature  $$\sigma_b:=(0; 2, \stackrel{b}{\ldots}, 2, 4,4,2^{n-1})$$where $b$ is an odd integer. We have detected a small mistake in this theorem since $b$ need not be odd; the following theorem amends the aforementioned mistake. 

\begin{theo} \label{corre}Let $n \geqslant 3$ be an integer and let $S$ be a compact Riemann surface with a group of automorphisms $G$ isomorphic to $Q(2^n).$ The following statements are equivalent.
\begin{enumerate} 
\item The action of $G$  on $S$ is a genus-zero action.
\item The signature of the action is $\sigma_b$ for some $b \geqslant 0.$ 
\end{enumerate}Moreover, for each $b \geqslant 0$ there exists a compact Riemann surface with a group of automorphisms isomorphic to $Q(2^n)$ acting on it with signature $\sigma_b;$ its genus is $2^{n-2}(b+1).$
\end{theo}

\subsection{Families of Riemann surfaces with $Q(2^n)$-action}\label{lugar}
 Hidalgo and Quispe in \cite{HQ} studied Riemann surfaces  with a triangle action (or, equivalently, regular dessins d’enfants in Grothendieck's terminology \cite{Grot}) of a group isomorphic to the so-called dicyclic group  $$G_m:=\langle x, y : x^{2m}, y^2x^{-m}, yxy^{-1}x\rangle \,\, \mbox{ where } \,\, m \geqslant 2.$$ 

For the special case of the  generalized quaternion group $Q(2^{n})=G_{2^{n-2}}$, they were able to prove that there exists, up to isomorphism, a unique Riemann surface, hereafter denoted by $X_n$, endowed with a triangle action of a group isomorphic to $Q(2^n).$ They also obtained that the signature of the action is  $(0; 4,4, 2^{n-1})$ (and therefore its genus is $2^{n-2}$), determined its full automorphism group, and proved that the surface is  represented by the  singular plane affine algebraic curve \begin{equation*}\label{lamp}Y^2 = X(X^{2^{n-1}} - 1).\end{equation*}

Observe that the hyperellipticity of $X_n$ coupled with Corollary \ref{isog-sup-gen} (even though this corollary was stated for $n \geqslant 4$, we shall see in \S\ref{TN} that it also holds for $n=3$) allow us to recover the fact proved in \cite[Theorem 2.1]{HQ} that each  group algebra decomposition of $JX_n$ with respect to $Q(2^n)$ is trivial.

\s

Let $\mathscr{M}_g$ denote the moduli space of compact Riemann surfaces of genus $g \geqslant 2$. It is known that $\mathscr{M}_g$ is a complex space of dimension $3(g-1)$ and that if $g \geqslant 4$ then $$\mbox{Sing}(\mathscr{M}_g)=\{[S]: S \mbox{ has nontrivial automorphisms}\}.$$

If a Riemann surface of genus $g$ admits a triangle action then it represents a {\it rigid} point in $\mbox{Sing}(\mathscr{M}_g)$; namely, it cannot be deformed together with its automorphisms. Hence, it is natural  to study  surfaces in $\mbox{Sing}(\mathscr{M}_g)$ that admit deformations keeping their automorphisms invariant. Among them, the simplest ones are the uniparametric families.

\s

A closed {\it family} $\mathcal{F}$ of  Riemann surfaces of genus $g$ is  a sublocus of $\mbox{Sing}(\mathscr{M}_g)$ consisting of all those surfaces that admit an action of a fixed group $K$ with a fixed signature $\sigma.$ Following \cite[Theorem 2.1]{b}, the interior of the family consists of finitely many {\it equisymmetric strata} that are in correspondence with the pairwise non-equivalent topologically actions of $K$ (see  \S\ref{preli3}). If $\sigma=(\gamma; k_1, \ldots, k_l)$  then the complex dimension of the family is given by $$\dim (\mathcal{F}) =3\gamma-3+l.$$

Observe that, in this terminology, Riemann surfaces with a triangle action correspond to  zero-dimensional families and among them, as already mentioned, the ones with a group of automorphisms isomorphic to $Q(2^n)$ have been classified in \cite{HQ}. Here we shall provide a  classification of  the  complex one-dimensional situation.

\begin{theo} \label{t-fam-sup} Let $n \geqslant 4$ be an integer. There exist precisely $n+1$ complex one-dimensional closed families of compact Riemann surfaces endowed with a group of automorphisms $G$ isomorphic to $Q(2^n).$ These families are:
\s
\begin{enumerate}
\item the family $\mathcal{F}_{n,0}$ consisting of all those compact Riemann surfaces of genus $2^{n-1}-1$ with $G$ acting on them with signature $(1; 2^{n-2}),$
\s
\item the family  $\mathcal{F}_{n,1}$ consisting of all those compact Riemann surfaces of genus $2^{n-1}+1$ with $G$ acting on them with signature  $(0; 4,4,4,4),$
\s
\item the family  $\mathcal{F}_{n,2}$ consisting of all those compact Riemann surfaces of genus $3 \cdot 2^{n-2}-1$ with $G$ acting on them with signature  $(0; 2^{n-1}, 2^{n-1}, 4,4),$ and
\s
\item  the family  $\mathscr{C}_{n,k}$, with $k \in \{2, \ldots, n-1\}$,  consisting of all those compact Riemann surfaces of genus $3 \cdot 2^{n-2}-2^{k-1}$ with $G$ acting on them with signature  $(0; 2^{n-1}, 2^{n-k}, 4,4).$ 
\end{enumerate}

\s

Furthermore, the families $\mathcal{F}_{n,0}, \mathcal{F}_{n,1}$ and $\mathscr{C}_{n,n-1}$ consist of only one equisymmetric stratum each, the family $\mathcal{F}_{n,2}$ consists of at most $2^{n-2}$ equisymmetric strata and the family $\mathscr{C}_{n,k}$ consists of at most $2^{n-k-1}$ equisymmetric strata, for each $k \in \{2, \ldots, n-2\}$.
\end{theo}

Given a compact Riemann surface with a group of automorphisms, a natural problem that arises is to decide whether or not the surface admits more automorphisms. This  is a challenging problem whose answer depends on the signature  of the action (as proved by Singerman in \cite{singerman2}) and on the geometry of a fundamental domain for the surface.

\s

For each $n \geqslant 4$ we consider the supergroups of $Q(2^n)$ of order $2^{n+1}$ given by
$$G_1=\langle x,y,z : x^{2^{n-1}}, z^2, y^2x^{2^{n-2}}, yxy^{-1}x, zxzx^{-1}, zyzy\rangle$$$$G_2=\langle x,y,z : x^{2^{n-1}}, z^2, y^2x^{2^{n-2}}, yxy^{-1}x, zxzx^{2^{n-2}+1}, zyzy\rangle.$$

 \begin{theo} \label{full-aut} Let $n \geqslant 4$ be an integer.
\s
\begin{enumerate}
\item If $S$ belongs to the family  $\mathcal{F}_{n,0}$ then the action of $G$ on  $S$ extends to an action of $G_1$ with signature $(0; 2,2,2,2^{n-1}).$ Furthermore, up to possibly finitely many exceptions, the full automorphism group of $S$ is isomorphic to $G_1.$
\s
\item If $S$ belongs to  the family $\mathcal{F}_{n,1}$ then  the action of $G$ on $S$ extends to an action of $G_1$ with signature $(0; 2,2,4,4).$ Furthermore, up to possibly finitely many exceptions, the full automorphism group of $S$ is isomorphic to $G_1.$
\s
\item If $S$ belongs to  the family $\mathcal{F}_{n,2}$ and the action of $G$ on $S$ extends to an action of a supergroup $G'$ of $G$ then $G'$ equals  either $G_1$ or $G_2$, and the signature of the action is $(0; 2,2,4,2^{n-1}).$ Furthermore, up to possibly finitely many exceptions, the full automorphism group of $S$ is isomorphic to either $G, G_1$ or $G_2.$  

\s
\item If $S$ belongs to  the family $\mathscr{C}_{n,k}$ then, up to finitely many exceptions, the full automorphism group of $S$ is isomorphic to $G,$ for each $k \in \{2, \ldots, n-1\}$.
\end{enumerate}
\end{theo}

\s

A consequence of the previous result together with Theorem \ref{corre} is that among the  previous families only $\mathscr{C}_{n,n-1}$ consists of Riemann surfaces with a genus-zero action of $Q(2^n).$ Hence, by Corollary \ref{isog-sup-gen}, the Jacobian varieties of the surfaces lying in the  remaining families admit a group algebra decomposition that is nontrivial. 

The following theorem provides such decompositions explicitly. 

\s

\begin{theo} \label{des-fam-0} Let $n \geqslant 4$ be an integer. 

\begin{enumerate}
\item If $S$ belongs to the family $\mathcal{F}_{n,0}$ then $JS_{N_1}, JS_{N_2}$ and $JS_{N_3}$ are elliptic curves, $$\dim \Pm(S \to S_Z)=2^{n-2}   \,\, \mbox{ and } \,\, \dim \Pm(S_{H_j} \to S_{H_{j+1}})=2^{n-j-2}$$for each $j \in \{2, \ldots, n-2\}.$ In particular $$JS \sim E \times \Pm(S \to S_Z) \times \Pi_{j=2}^{n-2} \Pm(S_{H_j} \to S_{H_{j+1}})^2$$where $E$ is an elliptic curve isogenous to $JS_G.$

\s

\item If $S$ belongs to the family $\mathcal{F}_{n,1}$ then $$\dim \Pm (S \to S_Z)=2^{n-1} \,\, \mbox{ and }\,\, \dim JS_{N_1}=1.$$ Moreover, $S_{N_2}, S_{N_3}$ and $S_{H_j}$ for each $j \in \{2, \ldots, n-2\}$  have genus zero, and$$JS \sim E \times \Pm(S \to S_{Z})$$where $E$ is an elliptic curve isogenous to $JS_{N_1}$. In particular 
$JS_{Z} \sim JS_{N_1}$.

\s

\item If  $S$ belongs to the family $\mathcal{F}_{n,2}$  then $S_{N_1}$ and $S_{N_2}$ have genus zero, $$\dim \Pm(S \to S_Z) = 2^{n-1} \,\, \mbox{ and }\,\, \dim \Pm(S_{H_j} \to S_{H_{j+1}})=2^{n-2-j}$$for each $j \in \{2, \ldots, n-2\}.$ Furthermore, if $n=4$ then $$JS \sim E_1 \times E_2^2 \times \Pm(S \to S_Z),$$whilst if $n \geqslant 5$ then $$JS \sim E_1 \times E_2^2 \times \Pm(S \to S_Z) \times \Pi_{j=2}^{n-3} \Pm(S_{H_j} \to S_{H_{j+1}})^2$$where $E_1$ and $E_2$ are elliptic curves isogenous to $JS_{N_3}$  and $JS_{{H}_{n-2}}$ respectively.

\s

\item If $S$ belongs to the family $\mathscr{C}_{n,k}$ and $k \in \{2, \ldots, n-2\}$ then $$ \dim \Pm(S \to S_Z) = 2^{n-1}\, \mbox{ and } 
\dim JS_{H_j} = 
 2^{n-j-1}-2^{k-2}$$if $2 \leqslant j\leqslant n-k$ and zero otherwise. In particular $S_{N_1}, S_{N_2}$ and $S_{N_3}$ have genus zero. Furthermore,   if $k=n-2$ then $$JS \sim \Pm(S \to S_Z) \times  JS_{H_2}^2,$$whilst  if $k \leqslant n-3$ then $$JS \sim \Pm(S \to S_Z) \times \Pi_{j=2}^{n-k-1} \Pm(S_{H_j} \to S_{H_{j+1}})^2 \times JS_{H_{n-k}}^2.$$
 \end{enumerate}
\end{theo}

For conciseness of exposition our results have been stated for $n \geqslant 4.$ The case $n=3$ is slightly different and analogous results hold;  we shall  discuss it briefly  in  \S\ref{TN}. At this point, we should anticipate that, with the notations of Theorem \ref{t-fam-sup}, the families $$\mathcal{F}_{3,0}, \,\, \mathcal{F}_{3,1}=\mathcal{F}_{3,2}\,\, \mbox{ and } \,\, \mathscr{C}_{3,2}$$ also exist, and that there is no other complex one-dimensional family of Riemann surfaces with quaternion group action.

\subsection{An algebraic description of $\mathscr{C}_{n,n-1}$} Among the families introduced in Theorem \ref{t-fam-sup}, only the family $\mathscr{C}_{n,n-1}$  consists of surfaces that are hyperelliptic. This property allows us to describe its members algebraically as follows.

\begin{prop} \label{malg}
Let $n \geqslant 3$ be an integer and let $S$ be a Riemann surface lying in the family $\mathscr{C}_{n,n-1}.$ Then $S$ is isomorphic to the normalization of the singular affine plane algebraic curve $$Y^2=X(X^{2^{n-2}}-1)(X^{2^{n-2}}-t)(X^{2^{n-1}}-t)$$for some $t \in \mathbb{C}-\{0,1\}$. In addition, with respect to this model, the group of automorphisms  of $S$ isomorphic to $Q(2^n)$ is generated by the transformations $$(X,Y) \mapsto (\xi^2 X, \xi Y) \,\, \mbox{ and } \,\, (X,Y) \mapsto (\tfrac{\lambda}{X}, \eta t \tfrac{Y}{X^{2^{n-1}+1}})$$where $\xi=\exppp(\tfrac{2 \pi i}{2^{n-1}}),$ and $\eta$ and $\lambda$ are chosen to satisfy $\eta^2=-\lambda$ and $\lambda^{2^{n-2}}=t.$ 
\end{prop}

Note that in the proposition above if we take $t = -1$ then we obtain the curve $$Y^2=X(X^{2^n}-1)$$mentioned in \S\ref{lugar} and denoted there by $X_{n+1}$. Thereby, as a direct consequence of the previous proposition, we derive the following interesting result.

\begin{coro} Let $n \geqslant 3$ be an integer. The family $\mathscr{C}_{n,n-1}$ contains 
 the unique  Riemann surface that admits a triangle action of the generalized quaternion group of order $2^{n+1}.$
\end{coro}

It is worth remarking that the exceptional Riemann surfaces uncovered by Theorem \ref{full-aut} correspond to those surfaces that have {\it extra automorphisms}; these surfaces lie in the closure of the families.  Observe that the previous corollary provides an explicit surface with this property for the   family $\mathscr{C}_{n,n-1}$, for each $n \geqslant 3.$

\subsection{Families of abelian varieties with $Q(8)$-action} It is well-known that the moduli space $\mathcal{A}_g$ of isomorphism classes of principally polarized abelian varieties of dimension $g \geqslant 3$ has a structure of complex  space of dimension $g(g+1)/2$ inherited by the projection $$\mathscr{H}_g \to \mathcal{A}_g:=\mathscr{H}_g/\Sp(2g, \mathbb{Z}),$$where $\mathscr{H}_g$ denotes the Siegel upper half-plane and $\Sp(2g, \mathbb{Z})$  the symplectic group. 

Moreover, the singular locus of $\mathcal{A}_g$ is given by
$$\mbox{Sing}(\mathcal{A}_g)=\{[A]: A \mbox{ has automorphisms different from}\pm id\}.$$

Here we construct three families of principally polarized abelian varieties of low dimension  whose members have a group of automorphisms isomorphic to the quaternion group. These families have the property that contain the three complex one-dimensional families of Jacobians with a quaternion group action.

\begin{theo}\label{ppavf1} There exists a
complex one-dimensional family $\mathscr{Y}_3$ of principally polarized abelian varieties of dimension three endowed with a group of  automorphisms isomorphic to $$\langle a,b,c : a^2,b^2, c^4, bcbc^3, acac^3, abac^2b\rangle =  (C_4 \times C_2) \rtimes C_2$$satisfying the following properties. 
\begin{enumerate}
\item The quaternion group of order eight acts on each member of   $\mathscr{Y}_3$ with analytic representation equivalent to $\chi_1 \oplus \Theta_1.$
\item The members of  $\mathscr{Y}_3$ decompose isogenously as the product of an elliptic curve and an abelian surface.
\item  $\mathscr{Y}_3$ contains the image under the Torelli map of the family $\mathcal{F}_{3,0}$.
\item $\mathscr{Y}_3$ is described, in terms of period matrices, as $$\left\{\left( \begin{smallmatrix}
1+i+\tfrac{1-i}{2}t & 1-it & \tfrac{1+i}{2}t \\
1-it & 1-(1+i)t & t \\
\tfrac{1+i}{2}t & t & i+\tfrac{i-1}{2}t
\end{smallmatrix} \right): t \in \mathbb{C}\right\}\cap \mathscr{H}_3$$where $i=\sqrt{-1}.$\end{enumerate}
\end{theo}

\begin{theo}\label{ppavf2} There exists a
complex two-dimensional family $\mathscr{Y}_5$ of principally polarized abelian varieties of dimension five endowed with a group of  automorphisms isomorphic to
$$\langle r,s,a,b : r^4,s^2, a^2, b^2, (sr)^2, arar^{-1}, (as)^2, 
brbr^{-1}, bsb(sra)^{-1}, bab(ar^2)^{-1}  \rangle = (\mathbb{D}_4 \times C_2) \rtimes C_2$$satisfying the following properties. 
\begin{enumerate}
\item The quaternion group of order eight  acts on each member of   $\mathscr{Y}_5$ with analytic representation equivalent to $\chi_2 \oplus 2\Theta_1.$
\item The members of  $\mathscr{Y}_5$ decompose isogenously as the product of an elliptic curve and an abelian fourfold.
\item  $\mathscr{Y}_5$ contains the image under the Torelli map of the family $\mathcal{F}_{3,1}=\mathcal{F}_{3,2}$.
\item $\mathscr{Y}_3$ is described, in terms of period matrices, as $$\left\{\left( \begin{smallmatrix}
\tfrac{3}{2}+4t_1-3t_2 & -\tfrac{7}{4}-2t_1+\tfrac{3}{2}t_2 & 0 & \tfrac{5}{4}-\tfrac{3}{2}t_2+2t_1 & \tfrac{7}{4}+\tfrac{5}{2}t_2-2t_1 \\
-\tfrac{7}{4}-2t_1+\tfrac{3}{2}t_2 & 1+t_1 & t_2 & -\tfrac{1}{2}+t_2-t_1 & \tfrac{1}{4}-\tfrac{1}{2}t_2+t_1\\
0 & t_2 & 2t_2 & \tfrac{1}{2} & t_2\\
\tfrac{5}{4}-\tfrac{3}{2}t_2+2t_1 & -\tfrac{1}{2}+t_2-t_1 & \tfrac{1}{2} & t_1 & -\tfrac{1}{4}+\tfrac{3}{2}t_2-t_1\\
\tfrac{7}{4}+\tfrac{5}{2}t_2-2t_1 & \tfrac{1}{4}-\tfrac{1}{2}t_2+t_1 & t_2 & -\tfrac{1}{4}+\tfrac{3}{2}t_2-t_1 & t_1
\end{smallmatrix} \right): t_1, t_2 \in \mathbb{C}\right\}\cap \mathscr{H}_5$$\end{enumerate}
\end{theo}

\s

\begin{theo}\label{ppavf3} There exists a
complex one-dimensional family $\mathscr{Y}_4$ of principally polarized abelian varieties of dimension four endowed with a group of  automorphisms $G$ isomorphic to the quaternion group of order eight satisfying the following properties. 
\begin{enumerate}
\item $G$ acts on each member of   $\mathscr{Y}_4$ with analytic representation  equivalent to $2 \Theta_1.$
\item The isotypical decomposition of each member of  $\mathscr{Y}_4$  with respect to $G$ is trivial.
\item  $\mathscr{Y}_4$ contains the image under the Torelli map of the family $\mathscr{C}_{3,2}.$
\end{enumerate}
\end{theo}

As proved in \cite{frediani}, the families $\mathcal{F}_{3,0}$ and $\mathscr{C}_{3,2}$ give rise to  {\it special subvarieties} in $\mathcal{A}_3$ and $\mathcal{A}_4$, and therefore their closures contain a dense set of points that are of CM type. 

\s

Consider the Riemann surface $X_4 \in \mathscr{C}_{3,2}$ and its Jacobian variety $JX_4 \in \mathscr{Y}_4.$ 
\begin{enumerate}

\item The full automorphism group of $X_4$ is isomorphic to the quasi-dihedral group $Q\mathbb{D}_{16}$ of order 32 and acts on $X_4$ with signature $(0; 2,4,16)$; see \cite[Section \S2.1.1]{HQ}.
\item $JX_4$ is of CM type; see, for instance,  \cite[Section 6.4]{W}.
\item $X_4$ is the Riemann surface with  full automorphism group of second largest order among the hyperelliptic surfaces of genus four; see the main result of \cite{Betal}.
\end{enumerate}

Even though we are not in position to describe the family $\mathscr{Y}_4$ as a whole in terms of period matrices as done for $\mathscr{Y}_3$ and $\mathscr{Y}_5$, we obtain the following result.

\begin{prop}\label{especial} A period matrix of $JX_4$ is given by$$\left( \begin{smallmatrix}
\tfrac{6}{7} - \tfrac{1}{28}\sqrt{2} + 23i \lambda & (10i\sqrt{2} - 27i)\lambda - \tfrac{3}{28}\sqrt{2} + \tfrac{1}{14} & (i\sqrt{2} - 5i)\lambda - \tfrac{3}{28}\sqrt{2} + \tfrac{1}{14} & (-13i\sqrt{2} + 19i)\lambda - \tfrac{1}{28}\sqrt{2} - \tfrac{1}{7}\\
(10i\sqrt{2} - 27i)\lambda - \tfrac{3}{28}\sqrt{2} +\tfrac{1}{14} & (-2i\sqrt{2} + 33i)\lambda + \tfrac{5}{28}\sqrt{2} + \tfrac{5}{7} & (3i\sqrt{2} - 15i)\lambda + \tfrac{5}{28} \sqrt{2} - \tfrac{2}{7}  & (i\sqrt{2} - 5i)\lambda - \tfrac{3}{28}\sqrt{2} + \tfrac{1}{14}\\
(i\sqrt{2} - 5i)\lambda - \tfrac{3}{28}\sqrt{2} + \tfrac{1}{14} & (3i\sqrt{2} - 15i)\lambda+ \tfrac{5}{28}\sqrt{2} - \tfrac{2}{7} & (-2i\sqrt{2} + 33i)\lambda + \tfrac{5}{28}\sqrt{2} + \tfrac{5}{7} & (10i\sqrt{2} - 27i)\lambda - \tfrac{3}{28}\sqrt{2} + \tfrac{1}{14}\\
(-13i\sqrt{2} + 19i)\lambda - \tfrac{1}{28}\sqrt{2} - \tfrac{1}{7} &  (i\sqrt{2} - 5i)\lambda - \tfrac{3}{28}\sqrt{2} + \tfrac{1}{14} & (10i\sqrt{2} - 27i)\lambda - 3/28\sqrt{2} + \tfrac{1}{14} &  \tfrac{6}{7} - \tfrac{1}{28}\sqrt{2} + 23i\lambda
\end{smallmatrix} \right)$$where $\lambda:=\tfrac{\sqrt{500 + 146\sqrt{2}}}{644}.$ 
\end{prop}

Section \S\ref{pre} is devoted to introducing some notations and to briefly reviewing the basic background. The case $n=3$ is  worked out in Section \S\ref{TN}. The proof of the results are given in the remaining sections.

\section{Preliminaries}\label{pre}

\subsection{Abelian varieties}\label{cri} Let $G$ be a finite group let $W$ be a rational irreducible representation of $G$. Then there is a complex irreducible representation $V$ of $G$ such that  $$W \otimes_{\mathbb{Q}} \mathbb{C}:=s_V (\oplus_{\sigma} V^{\sigma}),$$ where $s_V$ is the Schur index of $V$ and the sum is taken over the Galois group associated to the character field $K_V$ of $V$.  We say that $V$ is associated to $W$.

\s

Assume that $W_1, \ldots, W_r$ are the pairwise non-equivalent rational irreducible representations of $G.$ Let $V_l$ be a  complex irreducible representation of $G$ associated to $W_l$ and  set
$$e_l=\tfrac{d_{V_l}}{|G|} \Sigma_{g \in G} \mbox{tr}_{K_{V_l} |\mathbb{Q}} (\chi_{V_l}(g^{-1}))g \in \mathbb{Q}[G],$$where $d_{V_l}$ is the degree of $V_l$. As proved in \cite{LR}, the decomposition of $1$ as the sum \begin{equation}\label{noche1}1=e_1 + \cdots + e_r \,\, \mbox{ in }\,\, \mathbb{Q}[G]\end{equation}yields the isotypical decomposition of $A$ with respect to $G$$$JS \sim A_{1} \times \cdots \times A_{r} \,\, \mbox{ where}\,\, A_{l} := \Xi (\alpha_l e_l)(A)$$ and $\alpha_l \in \mathbb{Z}$ is chosen  to satisfy that $\alpha_l e_l \in {\mathbb Z}[G]$. In addition,  there are  (in general, non-uniquely determined by $W_l$) idempotents $f_{l1},\dots, f_{ln_l}$ such that \begin{equation}\label{noche2} e_l=f_{l1}+\dots +f_{ln_l} \, \, \mbox{ for each } \, 1 \leqslant l \leqslant r\end{equation}where  $n_l=d_{V_l}/s_{V_l}$. These idempotents provide $n_l$ pairwise isogenous subvarieties of $A.$ If we denote by  $B_l$  one of them for each $l,$ then \eqref{noche1} and \eqref{noche2} provide a group algebra decomposition of $A$ with respect to $G$
\begin{equation*} \label{eq:gadec}
A \sim B_{1}^{n_1} \times \cdots \times B_{r}^{n_r}.
\end{equation*}
See also \cite{CR} and \cite{RCR} for the special case of Jacobians.

Decompositions of abelian varieties with respect to special classes of groups can be found, for instance, in  \cite{d1},  \cite{nos},  \cite{PA}, \cite{jpaa}, \cite{pisa}, \cite{d2} and \cite{d4}.

\subsection{Riemann surfaces and Fuchsian groups}  The algebraic structure of a (co-compact) Fuchsian group $\Delta$ is determined by its {\it signature}; namely, the tuple\begin{equation*} \label{sig} \sigma(\Delta)=(\gamma; k_1, \ldots, k_s),\end{equation*}where  $\gamma$ is the genus of  $\mathbb{H}_{\Delta}$, with $\mathbb{H}$ denoting the upper half-plane, and $k_1, \ldots, k_s$ are the branch indices in the canonical projection $\mathbb{H} \to \mathbb{H}_{\Delta}.$ In this case, $\Delta$ has a presentation 
\begin{equation*}\label{prese}\langle \alpha_1, \ldots, \alpha_{\gamma}, \beta_1, \ldots, \beta_{\gamma}, x_1, \ldots , x_s : x_1^{k_1}=\cdots =x_s^{k_s}=\Pi_{i=1}^{\gamma}[\alpha_i, \beta_i] \Pi_{i=1}^s x_i=1\rangle,\end{equation*}where the brackets stand for the commutator.

Let $\Delta'$ be a group of automorphisms of $\mathbb{H}$ such that $\Delta \leqslant \Delta'$ of  finite index. Then $\Delta'$ is also Fuchsian and they are related by the so-called Riemann-Hurwitz formula $$ \mu(\Delta)/\mu(\Delta')= [\Delta' : \Delta],$$where $\mu(\Delta)=2\gamma-2 + \Sigma_{i=1}^s(1-1/{k_i}).$

\s

Let $S$ be a Riemann surface of genus $g \geqslant 2.$ By the  uniformization theorem, there is a Fuchsian group $\Gamma$ of signature $(g; -)$ such that $S \cong \mathbb{H}_{\Gamma}.$ Moreover, a finite group $G$ acts on $S$ if and only if there is a Fuchsian group $\Delta$  together with a group  epimorphism \begin{equation}\label{episs}\theta: \Delta \to G \, \, \mbox{ such that }  \, \, \mbox{ker}(\theta)=\Gamma.\end{equation}

This important result is known as Riemann's existence theorem; see, for example, \cite{Brou}, \cite{anita-ibero} and \cite{singerman}. Note that $G$ acts on $S$ with signature $\sigma(\Delta)$. We shall say that  the action is represented by the {\it surface-kernel epimorphism}  $\theta$ (hereafter {\it ske} for short). 

\s

Let $G'$ be a finite supergroup of $G.$ The action of $G$ on $S$ represented by the ske \eqref{episs} {\it extends} to an action of $G'$ on $S$ if there is a Fuchsian group $\Delta'$ together with a ske  $$\Theta: \Delta' \to G' \, \, \mbox{ in such a way that }  \, \, \Theta|_{\Delta}=\theta  \mbox{ and } \mbox{ker}(\theta)=\mbox{ker}(\Theta),$$ and the Teichm\"{u}ller spaces of $\Delta$ and $\Delta'$ have the same dimension. Maximal actions are those that cannot be extended. We refer to \cite{singerman2} and \cite{singerman} for more details.

\subsection{Equisymmetric stratification of $\mathscr{M}_g$} \label{preli3}Two actions $\epsilon_1, \epsilon_2$ of  $G$ on $S$  are  {\it topologically equivalent} if there exist $\omega \in \Aut(G)$ and $\varphi \in \mbox{Hom}^{+}(S)$  such that
\begin{equation*}\label{equivalentactions}
\epsilon_2(g) = \varphi \epsilon_1(\omega(g)) \varphi^{-1} \, \mbox{ for all } \, g\in G.
\end{equation*}

Note that each $\varphi$ as before induces an automorphism $\varphi^*$ of $\Delta$ where $\mathbb{H}_{\Delta} \cong S_G$. If $\mathscr{B}$ is the subgroup of $\mbox{Aut}(\Delta)$ consisting of them, then $\mbox{Aut}(G) \times \mathscr{B}$ acts on the set of skes defining actions of $G$ on $S$ with signature $\sigma(\Delta)$ by $$((\omega, \varphi^*), \theta) \mapsto \omega \circ \theta \circ (\varphi^*)^{-1}.$$  

Two skes $\theta_1, \theta_2 : \Delta \to G$ define topologically equivalent actions if and only if they belong to the same orbit of the aforementioned action; see, for example, \cite{Brou}. If $S_G \cong \mathbb{P}^1$ then $\mathscr{B}$ is known to be generated by the {\it braid transformations}  $\Phi_{i}$, for $1 \leqslant i  < s,$ defined by \begin{equation*} \label{braid} x_i \mapsto x_{i+1}, \hspace{0.3 cm}x_{i+1} \mapsto x_{i+1}^{-1}x_{i}x_{i+1} \hspace{0.3 cm} \mbox{ and }\hspace{0.3 cm} x_j \mapsto x_j \mbox{ when }j \neq i, i+1.\end{equation*}

Following \cite{b} (see also \cite{Harvey}), the singular locus of  $\mathscr{M}_g$ admits an  {\it equisymmetric} stratification where 
each equisymmetric stratum, if nonempty, corresponds to one topological class of maximal actions. More specifically:
$$\mbox{Sing}(\mathscr{M}_g)= \cup_{G \neq 1, \theta} \bar{\mathscr{M}}_g^{G, \theta}$$where the {\it equisymmetric stratum} ${\mathscr{M}}_g^{G, \theta}$ consists of surfaces of genus $g$ with full automorphism group isomorphic to $G$ such that the action is topologically equivalent to $\theta$. In addition, the  {\it closure} $\bar{\mathscr{M}}_g^{G, \theta}$ of  ${\mathscr{M}}_g^{G, \theta}$ is a closed irreducible algebraic subvariety of $\mathscr{M}_g$ and consists of surfaces  of genus $g$ with a group of automorphisms isomorphic to $G$ such that the action is  topologically equivalent to $\theta$.

\subsection{The moduli space of abelian varieties}\label{esp}   Each abelian variety $A=V/\Lambda$ admits a polarization;  namely, a non-degenerate real alternating form $E$ on $V$ such that $$E(\Lambda \times \Lambda) \subset \mathbb{Z} \,\,  \mbox{ and } \,\, E(iv, iw)=E(v,w) $$for all $v,w \in V$. If the elementary divisors of $E|_{\Lambda \times \Lambda}$ equal 1 then
$E$ is called {\it principal} and  $A$ is called a
 {\it principally polarized abelian variety} ({\it ppav} for short).  In this case,  there exists a basis for $\Lambda$ such that the matrix for $E|_{\Lambda \times \Lambda}$ with respect to it is given by \begin{equation}\label{simpl}
J = \left( \begin{smallmatrix}
0 & I_g \\
-I_g & 0
\end{smallmatrix} \right) \,\, \mbox{ where } \,\, g= \dim(A).
\end{equation}In addition, there exist a basis for $V$ with respect to which the period matrix for $A$ has the form $(I_g \, Z)$ where $Z  \in \mathscr{H}_g.$ If $(I_{g} \, Z_i)$ is the period matrix of $A_i$ then an isomorphism $A_1 \to A_2$ is given by invertible matrices \begin{equation}\label{ig}M \in \mbox{GL}(g, \mathbb{C}) \,\, \mbox{ and } \,\, R \in \mbox{GL}(2g, \mathbb{Z}) \,\, \mbox{ such that } \,\,  M(I_{g} \, Z_1)=(I_{g} \, Z_2)R.\end{equation}

As $R$ preserves the polarization \eqref{simpl}, it belongs to the  symplectic group and therefore  $$\mbox{Sp}(2g, \mathbb{Z}) \times \mathscr{H}_g \to \mathscr{H}_g \, \mbox{ given by } \,   (R=  \left( \begin{smallmatrix}
A & B \\
C & D
\end{smallmatrix} \right) , Z ) \mapsto R \cdot Z := (AZ+B)(CZ+D)^{-1}$$is an action 
whose orbits are in bijection with the isomorphism classes of ppavs. Hence$$\mathscr{H}_g \to \mathcal{A}_g:=\mathscr{H}_g/ \mbox{Sp}(2g, \mathbb{Z})$$is the moduli space of  of ppavs of dimension $g.$ See \cite{oort}.

\section{Proof of Theorem \ref{isog-sup-gen} and Proposition \ref{pp1}}

\subsection{Rational irreducible representations} \label{13} We recall here that  the complex irreducible representations of  $Q(2^n)$ are, up to equivalence, \begin{equation*}
\Theta_s  : x \mapsto 
 \mbox{diag} (\omega^s, \bar{\omega}^s) \, \mbox{ and }\, y \mapsto \left( \begin{smallmatrix}
0 & (-1)^s \\
1 & 0 \\
\end{smallmatrix} \right)
\end{equation*}where $\omega:=\mbox{exp}(2 \pi i / 2^{n-1})$ and $s \in  \{1, \ldots, 2^{n-2}-1\}$, and
$$\chi_1 : x, y \mapsto 1, \,\, \chi_2: x \mapsto 1, y \mapsto -1, \,\, \chi_3: x \mapsto -1, y \mapsto 1 \,\, \mbox{ and }\,\, \chi_4: x,y \mapsto -1.$$

For each $l \in \{1, \ldots, n-2\},$ the Galois group of the character field of $\Theta_{2^{l-1}}$ is $$\mathbb{G}_l =\mbox{Gal}(\mathbb{Q}(\omega^{2^{l-1}} + \bar{\omega}^{2^{l-1}})/\mathbb{Q}) = \mbox{Gal}(\mathbb{Q}( \cos(\tfrac{2 \pi}{2^{n-l}})/\mathbb{Q}) \cong C_{2^{n-l-2}}.$$

Observe that if $l \neq 1$ then $\Theta_{2^{l-1}}$ is not faithful,  induces a representation of the dihedral group of order $2^{n-1}$ and therefore its Schur index equals one. On the other hand, $\Theta_{1}$ is not real since its Schur-Frobenius indicator is equal to $-1.$ Thus, $\Theta_{1}$ cannot be defined over its character field, showing that its Schur index equals two. It then follows that the rational irreducible representations of $Q(2^n)$ are $\chi_1, \chi_2, \chi_3, \chi_4$ $$W_1:=2(\oplus_{\sigma \in \mathbb{G}_1} \Theta_{1}^{\sigma}) \,\mbox{ and } \, W_l := \oplus_{\sigma \in \mathbb{G}_l} \Theta_{2^{l-1}}^{\sigma}\mbox{ for } l \in \{2, \ldots, {n-2}\}.$$

\subsection{Proof of Theorem \ref{isog-sup-gen}}

Let $A$ be an abelian variety with a group of automorphisms $G$ isomorphic to $Q(2^n)$. Then the aforementioned facts concerning the rational irreducible representations of $Q(2^n)$ together with the  results mentioned in \S\ref{cri} show that the isotypical decomposition of $A$ with respect to $G$ is given by \begin{equation*} A \sim A_{G} \times A_2 \times A_3  \times A_4 \times A_{W_1} \times A_{W_2} \times \cdots \times A_{W_{n-2}}\end{equation*}where $A_i$ is the abelian subvariety of $A$ associated to $\chi_i$ and $A_{W_l}$ is the one associated to $W_l.$ Besides,  the group algebra decompositions of $A$ with respect to $G$ have the form \begin{equation}\label{isog-gen} A \sim A_{G} \times B_2 \times B_3  \times B_4 \times B_{W_1} \times B_{W_2}^2 \times \cdots \times B_{W_{n-2}}^2\end{equation}where $B_i=A_i$, $B_{W_1} = A_{W_1}$  and $B_{W_l}$ is a subvariety of $A$ such that $B_{W_l}^2 \sim A_{W_l}$ for $l \geqslant 2$.

\s

Let $K$ be a subgroup of $G$. In the sequel, we denote by $\rho_K$ the representation of $G$ induced by the trivial one of $K.$ Note that $\rho_G$ and $\rho_{1_G}$ are the trivial and  regular representation, respectively. We recall that,  by Frobenius Reciprocity theorem (see, for instance, \cite[Chapter  7]{Serre}), if $V$ is a complex representation of $G$ then $$\langle V, \rho_K  \rangle_G=\dim V^{K}$$ where $V^K$ denotes  the vector subspace of $V$ of elements fixed under $K$ and the brackets stand for the inner product of characters of $G.$ It then follows that $$\rho_K \simeq \oplus_{i=1}^4 (\dim \chi_i^{K}) \chi_i \oplus \oplus_{s=1}^{2^{n-2}-1}(\dim \Theta_s^K)\Theta_s.$$

In order to avoid redundant computations, we also mention that $\dim V^K$ is invariant under the action of the Galois group associated to the character field of $V$.

\s

Observe that the dimension of the fixed subspace of $\chi_i$ under the action of $H_j$ equals 1 if  $i= 1$ and $i=3$, and equals zero otherwise, for each $j.$ In addition, the dimension of the fixed subspace of $\Theta_1$ under the action of $H_j$ equals zero, for all $j.$ Similarly, for $l \geqslant 2$ 
 \begin{displaymath}
\Theta_{2^{l-1}}(x^{2^{n-j}})=   \mbox{diag}(
\omega^{2^{n-j+l-1}}, \bar{\omega}^{2^{n-j+l-1}})  \end{displaymath}and, consequently, the dimension of the fixed subspace of $\Theta_{2^{l-1}}$ under the action of $H_j$ is equal to 1 if $l \geqslant j$ and zero otherwise. All the above implies that  $$\rho_{H_j} \simeq \chi_1 \oplus \chi_3 \oplus \oplus_{l \geqslant j} ( \oplus_{\sigma \in G_{l}}\Theta_{2^{l-1}}^{\sigma}) = \chi_1 \oplus \chi_3 \oplus (\oplus_{l \geqslant j} W_l)$$and therefore $$\rho_{H_{j}} \simeq \rho_{H_{j+1}} \oplus W_j \,\, \mbox{ for each } \,\, j \in \{2, \ldots, n-2\}.$$ We now apply \cite[Theorem 6.2]{CR-H} to conclude that  there exists an isogeny\begin{equation} \label{i1}B_{W_j} \to \Pm(A_{H_j} / A_{H_{j+1}})\mbox{ for each } \,\, j \in \{2, \ldots, n-2\}.\end{equation}

By proceeding analogously, it can be seen that $$\rho_{1_G} \simeq \rho_{Z} \oplus W_1, \,\, \rho_{N_1} \simeq \rho_{G} \oplus \chi_2, \,\,  \rho_{N_2}\simeq \rho_{G}  \oplus \chi_3  \, \mbox{ and } \, \rho_{N_3} \simeq \rho_G \oplus \chi_4.$$Hence, there exist isogenies \begin{equation} \label{i2}B_{W_1} \to \mbox{Prym}(A/ A_{Z})\, \mbox{ and } \, B_{l+1} \to \Pm(A_{N_l} / A_{G})\end{equation}for each $l \in \{1,2,3\}$,  and the desired isogeny  $$A \sim A_{G} \times \Pi_{i=1}^3   \Pm(A_{N_i} /A_{G}) \times \Pm(A/ A_{Z}) \times \Pi_{j=2}^{n-2} \Pm(A_{H_j} / A_{H_{j+1}})^2$$is obtained by replacing \eqref{i1} and \eqref{i2} in \eqref{isog-gen}. 

\s

To obtain the remaining isogeny decompositions, we proceed analogously to see  $$\rho_{\tilde{H}_{j}} \simeq \rho_{\tilde{H}_{j+1}} \oplus W_j \,\, \mbox{ and } \,\, \rho_{K_i} \simeq  \rho_{K_{i+1}} \oplus 2W_i$$ for each $i,j \in \{2, \ldots, n-2\},$ and therefore $$B_{W_j} \sim\Pm(A_{\tilde{H}_j} / A_{\tilde{H}_{j+1}}) \,\, \mbox{ and } \,\,   A_{W_i} = \Pm(A_{K_i}/A_{K_{i+1}})\sim B_{W_i}^2.$$
\subsection{Proof of Proposition \ref{pp1}} Let $n \geqslant 4$ be an integer and let $A$ be an abelian variety with a group of automorphisms $G$ isomorphic to  $Q(2^n).$ Let $s \in \{1, \ldots, 2^{n-2}-1\}.$ Since the character of $\Theta_s$ is real, the complex-conjugate representation $\bar{\Theta}_s$ is equivalent to $\Theta_s$. If $\rho := \rho_a \oplus \bar{\rho}_a$  and $\sigma$ belongs to the Galois group of the character field of $\Theta_s$ then $$\langle \rho_a, \Theta_s\rangle_G =\tfrac{1}{2}\langle \rho, \Theta_s\rangle_G=\tfrac{1}{2}\langle \rho, \Theta_s^{\sigma}\rangle_G=\langle \rho_a, \Theta_s^{\sigma}\rangle_G,$$where the second equality follows from the fact that $\rho$ is rational. In other words, the multiplicity of the irreducible summands of $\rho_a$ is invariant under Galois-conjugation. Thus, having said that, to prove the first statement of the proposition it suffices to notice that the Galois orbits of $\Theta_1$ and of $\Theta_{2^{j-1}}$, for each $j \in \{2, \ldots, n-2\},$ are $$\{\Theta_1, \Theta_3, \dots, \Theta_{2^{n-2}-1}\} \,\, \mbox{ and }\,\, \{\Theta_{t_j^k2^{j-1}} : 0 \leqslant k < 2^{n-j-2}\}$$respectively,  where $t_j$ is as in the statement of the proposition. Note that, with the notations of \S\ref{13}, the analytic representation of $G$ can be rewritten as
$$\rho_a \simeq \oplus_{j=1}^4 a_j \chi_j \oplus b_1 (\oplus_{\sigma \in G_1} \Theta_{1}^{\sigma}) \oplus \oplus_{l=2}^{n-2} b_{2^{l-1}} W_l.
$$

To prove the second statement we recall that, following \cite[Corollary 5.4]{CR-H},  $$\dim A_K=\langle \rho_a, \rho_K\rangle_G=\dim \rho_a^{K}$$for each subgroup $K$ of $G.$ It follows that $\dim A_G=\langle \rho_a, \rho_G\rangle_G = \langle \rho_a, \chi_1 \rangle_G =a_1$ and $$\dim \Pm(A_{N_l}/A_G)=\dim A_{N_l}-a_1=\langle \rho_a, \rho_{N_l}\rangle_G-a_1=(a_{1}+a_{l+1})-a_1=a_{l+1}.$$
Similarly we have that $$\dim \Pm(A/A_Z)=\langle \rho_a, \rho_{1_G}\rangle_G-\langle \rho_a, \rho_Z\rangle_G=2(b_1+b_3+\cdots+b_{2^{n-2}-1})=2^{n-2}b_1,$$where the last equality follows from the fact that $b_s = b_1$ for each $s$ odd.

\s

Finally, for each $j \in \{2, \ldots, n-2\}$ we have that $$\dim \Pm(A_{H_j}/A_{H_{j+1}})=\langle \rho_a, \rho_{H_j}\rangle_G- \langle \rho_a, \rho_{H_{j+1}}\rangle_G=\Sigma_{k=0}^{2^{n-j-2}-1} b_{t_j^k2^{j-1}}=2^{n-j-2}b_{2^{j-1}}$$where the last equality follows from the fact that $b_{t_j^k2^{j-1}}=b_{2^{j-1}}$ for each $k.$ 

\section{Proof of Theorem \ref{equiv}}Let $n \geqslant 4$ be an integer and let $A$ be an abelian variety with a group of automorphisms $G$ isomorphic to  $Q(2^n).$ By Proposition \ref{isog-sup-gen},  the analytic representation of the action of $G$ on $A$ is given by  $$\rho_a \simeq \oplus_{j=1}^4 a_j \chi_j \oplus b_1 (\oplus_{\sigma \in G_1} \Theta_{1}^{\sigma}) \oplus \oplus_{l=2}^{n-2} b_{2^{l-1}} W_l.
$$ for some  non-negative integers $a_1, \ldots, a_4, b_1, b_{2}, b_4, \ldots, b_{2^{n-3}}.$

\s

(1) $\implies$ (2). Observe that if the dimension of $A_Z$ is positive then the dimension of $A$ is strictly greater than the dimension of $\Pm(A/A_Z)$ and therefore the group algebra decomposition of $A$ with respect to $G$ \begin{equation*}\label{eend}A \sim A_{G} \times \Pi_{i=1}^3   \Pm(A_{N_i} /A_{G}) \times \Pm(A/ A_{Z}) \times \Pi_{j=2}^{n-2} \Pm(A_{H_j} / A_{H_{j+1}})^2\end{equation*}is nontrivial.  

\s

(2) $\implies$ (3). If $K$ is any nontrivial subgroup of $G$ then $Z \leqslant K$ and therefore $A_K \subseteq A_Z.$ It then follows that if the dimension of $A_Z$ equals zero then the same conclusion holds for any nontrivial subgroup $K$ of $G.$


\s

(3) $\implies$ (4). It is a direct consequence of Proposition \ref{pp1}.

\s

(4) $\implies$ (5). Observe that if  $s$ is odd then $1$ is not an eigenvalue of the matrix $\Theta_s(g)$ for all nontrivial $g \in G.$ The same conclusion is valid for the representation $$ b_1( \Theta_1 \oplus  \Theta_3 \oplus \cdots \oplus  \Theta_{2^{n-2}-1})$$ for each non-negative integer $b_1$, which is assumed to be equivalent  to $\rho_a$. 

\s

(5) $\implies$ (1).  If we assume $\rho_a$ to be fixed point free then (\cite[Corollary 5.4]{CR-H}) $$\dim A_K=\langle \rho_a, \rho_K\rangle_G=\dim \rho_a^{K}=0$$for each nontrivial subgroup $K$ of $G,$ and therefore each group algebra decomposition of $A$ with respect $G$  becomes trivial.

\section{Proof of Theorem \ref{corre}}

Let $n \geqslant 3$ be an integer and let $S$ be a compact Riemann surface of genus $g \geqslant 2$ with a group of automorphisms $G$ isomorphic to $Q(2^n).$ Let $\pi$ denote the regular branched covering map given by the action of $G$ on $S.$ Observe that the signature of the action of $G$ on $S$ must have the following form: $$(\gamma; 4, \stackrel{a}{\ldots},4, 2, \stackrel{b}{\ldots}, 2, 2^{n-1}, \stackrel{c_1}{\ldots}, 2^{n-1}, \ldots, 2^{n-k}, \stackrel{c_k}{\ldots}, 2^{n-k}, \ldots, 2^{2}, \stackrel{c_{n-2}}{\ldots}, 2^2)$$where $a,b,c_1, \ldots, c_{n-2}$ are non-negative integers. Here, the branch points of $\pi$ associated to the $a$ branch values of $\pi$ marked with 4 have isotropy group conjugate to $\langle y \rangle$ or $\langle xy \rangle,$ whereas the ones associated to the the $c_{n-2}$ branch values of $\pi$ marked with 4 have isotropy group generated by $\langle x^{2^{n-3}}\rangle.$

\s

Assume the action of $G$ to be of genus-zero. Then $\gamma=0$ and, by Riemann's existence theorem, $a$ is even and different from zero. The Riemann-Hurwitz formula says that  \begin{equation}\label{ee1}2(g-1)=2^n(-2+\tfrac{3}{4}a+\tfrac{1}{2}b+\Sigma_{k=1}^{n-2}c_k(1-\tfrac{1}{2^{n-k}})).\end{equation}

We denote by $\pi_Z$ the regular branched covering map given by the action of $Z$ on $S.$ Observe that each branch value of $\pi$ marked with 4 yields $2^{n-2}$ branch values of $\pi_Z,$ each branch value of $\pi$ marked with 2 yields $2^{n-1}$ branch values of $\pi_Z,$ and each branch value of $\pi$ marked with $2^{n-k}$ yields $2^{k}$ branch values of $\pi_Z$ for each $k \in \{1, \ldots, n-2\}.$

It then follows that the signature of the action of $Z$ on $S$ is $$(0; 2, \stackrel{d}{\ldots},2) \,\, \mbox{ where } \, d=a2^{n-2}+b2^{n-1}+\Sigma_{k=1}^{n-2}c_k2^k$$and therefore we obtain that \begin{equation}\label{ee2}2(g-1)=-4+a2^{n-2}+b2^{n-1}+\Sigma_{k=1}^{n-2}c_k2^k\end{equation}

The equations \eqref{ee1} and \eqref{ee2} show that $$-2+\tfrac{3}{4}a+\tfrac{1}{2}b+\Sigma_{k=1}^{n-2}c_k(1-\tfrac{1}{2^{n-k}})=-2^{2-n}+\tfrac{1}{4}a+\tfrac{1}{2}b+\Sigma_{k=1}^{n-2}c_k2^{k-n}$$and therefore there is no restriction on $b.$ The previous equality can be rewritten as $$\Sigma_{k=1}^{n-2}c_k(1-\tfrac{1}{2^{n-k-1}})+\tfrac{1}{2}a=2-\tfrac{4}{2^n}.$$Observe that if $a \geqslant 4$ then $\Sigma_{k=1}^{n-2}c_k(1-\tfrac{1}{2^{n-k-1}}) <0$; a contradiction. Hence, $a=2$ and $$\Sigma_{k=1}^{n-2}c_k(1-\tfrac{1}{2^{n-k-1}})=\Sigma_{k=2}^{n-2}c_k(1-\tfrac{1}{2^{n-k-1}})+c_1(1-\tfrac{4}{2^{n}})=1-\tfrac{4}{2^n}.$$Thus, $c_1=1$ and $c_k =0$ for $k \neq 1,$ showing that the signature of the action is  $$\sigma_b=(0; 2, \stackrel{b}{\ldots},2,4,4,2^{n-1})$$for some $b \geqslant 0,$ as desired. Conversely, if the signature of the action of $G$ is $\sigma_b$ then the signature of the action of $Z$ on $S$ is of the form $$(\gamma; 2, \stackrel{e}{\ldots}, 2) \,\, \mbox{ where } e=2^{n-1}+b2^{n-1}+2,$$for some $\gamma \geqslant 0.$ By proceeding analogously as already done, one sees that  $$2^n(\tfrac{1}{2}+\tfrac{1}{2}b-\tfrac{1}{2^{n-1}})=4\gamma-4+e$$ and then $\gamma=0.$ It follows from Corollary \ref{isog-sup-gen}  that the action of $G$ is a genus-zero action.

\s

It only remains to prove that for each $b \geqslant 0$ there exists a  Riemann surface $S$ on which  $G$ acts with signature $\sigma_b.$ Let $\Delta$ be a Fuchsian group of signature $\sigma_b$ presented as$$\Delta=\langle x_1, \ldots, x_b, y_1, y_2, z_1 : x_1^2=\cdots=x_b^2=y_1^4=y_2^4=z_1^{2^{n-1}}= x_1 \cdots  x_by_1y_2z_1=1 \rangle.$$ The homomorphisms $\eta: \Delta \to Q(2^n)$ given by $$x_1 \mapsto y^2, \ldots, x_b \mapsto y^2 , y_1 \mapsto y, y_2 \mapsto y^{-1}x^{-1}, z_1 \mapsto x \,\, \mbox{ if } b \mbox{ is even} $$ $$x_1 \mapsto y^2, \ldots, x_b \mapsto y^2 , y_1 \mapsto x^{-1}y, y_2 \mapsto y, z_1 \mapsto x \,\, \mbox{ if } b \mbox{ is odd} $$are two skes. Hence, we can take $S:= \mathbb{H}/\mbox{ker}(\eta)$ and the result follows from Riemann's existence theorem. 

%

\section{Proof of Theorem \ref{t-fam-sup}}

Let $\Delta_0$ be a Fuchsian group of signature $\sigma_0=(0; k_1, k_2, k_3, k_4)$  presented as  \begin{equation}\label{dom}\Delta_0=\langle x_1, x_2, x_3, x_4 : x_1^{k_1}= x_2^{k_2}= x_3^{k_3}= x_4^{k_4}=x_1x_2x_3x_4=1 \rangle\end{equation}and let $\theta: \Delta_0 \to G=Q(2^n)$ be a ske representing an action of $G$ on a  Riemann surface $S$  with signature $\sigma_0$. Write $g_j=\theta(x_j)$ for each $j$ and identify $\theta$ with the $4$-tuple$$\theta=(g_1, g_2, g_3, g_4).$$

The surjectivity of $\theta$ shows that the number $t$ of elements $g_j$'s that lie in $\langle x \rangle$ is at most two. In addition, as the product of three elements of $G-\langle x \rangle$ does not belong to $\langle x \rangle$, we see that $t$ cannot equal 1. The situation naturally splits into three cases.

\s

{\bf Case A.} Assume $t=0.$ Then $k_j=4$ for each $j$ and the surfaces $S$ with this type of action belong to the family $\mathcal{F}_{n,1}$. If we write $g_j=x^{n_j}y$  then, after considering  suitable braid transformations and the automorphism $(x,y) \mapsto (x,xy)$ of $G$, one of the following cases occurs: each $n_j$ is even; 
 $n_1, n_2, n_3$ are even and $n_4$ is odd;  or
 $n_2, n_3$ are even and $n_1, n_4$ are odd. Observe that the first case yields  a non-surjective epimorphism $\theta,$ whereas the second case gives rise to a map which sends $x_1x_2x_3x_4$ to a nontrivial element of $G.$ Thus, again after applying an appropriate automorphism of $G$, the ske $\theta$ is equivalent to \begin{equation}\label{digi}\theta_p:=(xy,y,x^py,x^{p+1}y) \, \mbox{ where } p \in \mathbb{Z}_{2^{n-1}} \mbox{ is even.}\end{equation} The fact that $\Phi_3^2$ identifies $\theta_p$ with $\theta_{p+2}$ allows us to conclude that the family $\mathcal{F}_{n,1}$ consists of  only one equisymmetric stratum (represented by $\theta_0$), as desired.  If $g$ is the genus of $S$ then  the Riemann-Hurwitz formula reads $$2(g-1)=2^n[-2+4(1-\tfrac{1}{4})] \,\, \mbox{ showing that } \,\, g=2^{n-1}+1.$$

\s

From now on we assume $t=2.$ After considering suitable braid transformations, we can assume that $k_3=k_4=4$ and that $\theta$ is equivalent to the ske \begin{equation*}\label{abske}(x^{\alpha}, x^{\beta}, g_3,g_4) \mbox{ where }\alpha, \beta \in \mathbb{Z}_{2^{n-1}}-\{0\}.\end{equation*} 

\s

{\bf Case B.} Assume $g_3$ and $g_4$ to be conjugate. Up to  $(x,y) \mapsto (x,xy)$ and after applying an inner automorphism of $G$, we can assume $g_3=y$ and $g_4=x^{p}y$ where $p \in \mathbb{Z}_{2^{n-1}}$ is even. Observe that if both $\alpha$ and $\beta$ were even then $\theta$ would not be surjective. Thus, after applying the braid transformation $\Phi_1$ if necessary, we can assume $\alpha=1$ and therefore $k_1=2^{n-1}.$ The fact that $x_1x_2x_3x_4=1$ says that $\beta$ equals $p-1+2^{n-2}$ in $\mathbb{Z}_{2^{n-1}}$ and, consequently, $\theta$ is equivalent to \begin{equation} \label{luz}\theta_{p}:=(x, x^{p-1+2^{n-2} }, y,x^py) \mbox{ where } p \in \mathbb{Z}_{2^{n-1}} \mbox{ is even.}\end{equation} 

The fact that $p-1+2^{n-2}$ is odd shows that $k_2=2^{n-1}$ and therefore the surfaces $S$ with this type of action form the family $\mathcal{F}_{n,2}.$ Note that this family consists of at most $2^{n-2}$ equisymmetric strata, according to the possible choices of $p$ in \eqref{luz}. If $g$ is the genus of $S$ then  the Riemann-Hurwitz formula implies that $g=3 \cdot 2^{n-2}-1.$

\s

{\bf Case C.} Assume $g_3$ and $g_4$ to be non-conjugate.  Up to  $(x,y) \mapsto (x,xy)$ and after applying an inner automorphism of $G$ we can assume $g_3=y$ and $g_4=xy.$ The fact that $x_1x_2x_3x_4=1$ implies that $\beta$ equals $1+2^{n-2}-\alpha$ in $\mathbb{Z}_{2^{n-1}}$ and $\theta$ is equivalent to \begin{equation*}\label{emo}\theta_{\alpha}:=(x^{\alpha}, x^{1+2^{n-2}-\alpha}, y,xy) \mbox{ where } \alpha \in \mathbb{Z}_{2^{n-1}}-\{0\} \mbox{ and } \alpha \neq 1+2^{n-2}.\end{equation*}

Since $\Phi_1$ identifies $\theta_{\alpha}$ with $\theta_{1+2^{n-2}-\alpha}$, we can assume $\alpha$ odd and therefore $k_1=2^{n-1}$. Observe that $1+2^{n-2}-\alpha$ can equal each one of the even integers between $2$ and $2^{n-2},$ showing that $k_2$ can equal $2^{n-k}$ for each $k \in \{2, \ldots, n-1\}.$ Thus, once $k$ is fixed, the surfaces $S$ with this type of action form the family $\mathscr{C}_{n,k}$ for each $k$ as before. Besides, the family $\mathscr{C}_{n,k}$ consists of at most  $2^{n-k-1}$ equisymmetric strata, because they are in correspondence with the generators of the cyclic subgroup of $\langle x \rangle$ of order ${2^{n-k}}.$  The computation of the genus of each $S$ lying in $\mathscr{C}_{n,k}$ follows from the Riemann-Hurwitz formula.

\s

Finally, we consider the case in which the quotient by the action of $G$ has genus one. Let $\Delta_1$ be a Fuchsian group of signature $\sigma_1=(1; k_1)$  presented as \begin{equation} \label{sab}\Delta_1=\langle \alpha, \beta, \gamma: [\alpha,\beta]\gamma=\gamma^{k_1}=1 \rangle\end{equation}and let $\theta: \Delta_1 \to G=Q(2^n)$ be a ske representing an action of $G$ on a compact Riemann surface $S$  with signature $\sigma_1$. We write $a=\theta(\alpha), \, b=\theta(\beta)$ and $g=\theta(\gamma)$ and identify $\theta$ with $\theta=(a,b,g).$ We recall that, following \cite[Proposition 2.5]{Brou}, the maps  $$(\alpha, \beta, \gamma) \mapsto (\alpha, \beta\alpha, \gamma) \,\, \mbox{ and } \,\,(\alpha, \beta, \gamma) \mapsto (\alpha\beta, \beta,\gamma)$$ identify actions that are equivalent. The previous fact together with the surjectivity of $\theta$ allow us to assume $a$ and $b$ to belong to $G-\langle x \rangle.$ If $a$ and $b$ were conjugate then we might assume $a=y$ and $b=x^{p}y$ where $p$ is even. But the fact that the commutator subgroup of $G$ equals $\langle x^2 \rangle \cong C_{2^{n-2}}$  would imply that $\theta$ is not surjective; a contradiction. Thus, up to equivalence, $a=y, b=xy$ and therefore $g=x^2$ and $k_1=2^{n-2}.$ Hence, the surfaces $S$ with this type of action form the family $\mathcal{F}_{n,0}$ and the computation of their genera follows from the Riemann-Hurwitz formula.

\section{Proof of Theorem \ref{full-aut}}

Let $S$ belong to the family $\mathcal{F}_{n,0}.$   The action of $G$ on $S$ is represented by  the ske \begin{equation}\label{back} \Delta_1 \to G \,\, \mbox{ given by }\,\, (\alpha, \beta, \gamma) \mapsto (y, xy,x^2),\end{equation}where $\Delta_1$ is a Fuchsian group of signature $(1; 2^{n-2})$ as in \eqref{sab}. Let $$\Delta_1'=\langle y_1, y_2, y_3, y_4: y_1^2=y_2^2=y_3^2=y_4^{2^{n-1}}=y_1y_2y_3y_4=1 \rangle$$be a Fuchsian group of signature $(0; 2,2,2,2^{n-1})$ and consider the ske $\theta': \Delta_1' \to G_1$ $$(y_1, y_2, y_3, y_4) \mapsto (z,yz,xyz,xz).$$Observe that the elements $$\alpha':=y_1y_2, \,\, \beta':= y_4y_2 \,\, \mbox{ and } \gamma':=y_4^2$$generate a subgroup of $\Delta_1'$  isomorphic to $\Delta_1.$ The restriction of $\theta'$ to $\Delta_1$ is given by $$(\alpha', \beta', \gamma') \mapsto (y^{-1}, xy^{-1},x^2).$$

Note that the image of the restriction is $\langle x, y \rangle \cong G.$ Moreover, after considering the automorphism of $G$ given by $(x, y) \mapsto (x, y^{-1})$, the restriction  agrees with \eqref{back}. Hence the action of $G$ on $S$ extends to an action of $G_1$ with signature $(0; 2,2,2,2^{n-1})$, as claimed. If, in addition, $S$ belongs to the interior of $\mathcal{F}_{n,0}$ (namely, for all the members of the family up to finitely many exceptional surfaces) then  $\mbox{Aut}(S) \cong G'$ due to the fact that  the signature $(0; 2,2,2,2^{n-1})$ is maximal; see  \cite[Theorem 1]{singerman2}.

\s

Let $S$ belong to the family $\mathcal{F}_{n,1}.$  The action of $G$ on $S$ is represented by the ske\begin{equation}\label{green} \Delta_0\to G \,\, \mbox{ given by } (x_1, x_2, x_3, x_4) \mapsto (xy,y,y^{-1},xy^{-1})\end{equation} (we have chosen $p=2^{n-2}$ in \eqref{digi}) where $\Delta_0$ is a Fuchsian group of signature $(0; 4,4,4,4)$ as in \eqref{dom}. Let $$\Delta_0'=\langle y_1, y_2, y_3, y_4: y_1^2=y_2^2=y_3^4=y_4^4=y_1y_2y_3y_4=1 \rangle$$be a Fuchsian group of signature $(0; 2,2,4,4)$ and consider the ske $\theta': \Delta_0' \to G_1$ $$(y_1, y_2, y_3, y_4) \mapsto (zy,xyz,xy,y).$$Observe that the elements $${x}_1':=y_3, \,\,{x}_2':= y_4, \,\, {x}_3':=y_1y_4y_1 \,\, \mbox{ and } {x}_4':=y_2y_3y_2$$generate a subgroup of $\Delta_0'$  isomorphic to $\Delta_0.$ The restriction of $\theta'$ to $\Delta_0$ is given by $$({x}_1',{x}_2',{x}_3',{x}_4' ) \mapsto (xy,y,y^{-1},xy^{-1}),$$that is, the restriction agrees with \eqref{green}. This proves the first assertion. 

To prove the second assertion we only have to show that the action of $G_1$ on $S$ does not extend to any supergroup of $G_1$. We proceed by contradiction. If the action extends then, by \cite[Theorem 1]{singerman2}, there must exist a group $G'$ of order $2^{n+2}$ such that $$G \leqslant G_1 \leqslant G'$$ and that acts on $S$ with signature $(0; 2,2,2,4).$ It then follows that $G'$ is generated by  three involutions whose product has order four; we denote these involutions by $a,b$ and $c.$ As $\langle x \rangle$ is a characteristic subgroup of $G_1$ and $G_1$ is a normal subgroup of $G'$ we obtain that $a,b$ and $c$ define automorphisms of $\langle x \rangle$ and, consequently, $a,b$ and $c$ must commute. Hence, the product $abc$ has order two; a contradiction.

\s

Let $S$ belong to the family $\mathcal{F}_{n,2}.$ The action of $G$ on $S$ is represented by  the ske \begin{equation} \label{luz2} \theta_p: \Delta_0 \to G \,\, \mbox{ given by }(x_1, x_2, x_3, x_4) \mapsto(x, x^{p-1+2^{n-2} }, y,x^py)\end{equation}for some $p$ even, where $\Delta_0$ is of signature $(0; 2^{n-1}, 2^{n-1},4,4)$  as in \eqref{dom}.  Let $$\Delta_0''=\langle y_1, y_2, y_3, y_4: y_1^2=y_2^2=y_3^4=y_4^{2^{n-1}}=y_1y_2y_3y_4=1 \rangle$$be a Fuchsian group of signature $(0; 2,2,4,2^{n-1})$ and consider the skes $\theta_1'': \Delta_0'' \to G_1$ and $\theta_2'': \Delta_0'' \to G_2$ both defined by the rule$$(y_1, y_2,y_3,y_4) \mapsto (z,  zy,  xy^{-1},  x).$$The elements $${x}_1'':=y_2y_4y_2, \,\, {x}_2'':= y_4, \,\, {x}_3'':=(y_1y_2)y_3(y_1y_2)^{-1} \,\, \mbox{ and } {x}_4'':=y_2y_3y_2$$generate a subgroup of $\Delta_0''$  isomorphic to $\Delta_0.$ The restrictions of $\theta_1''$ and $\theta_2''$ to $\Delta_0$ are  \begin{equation} \label{resf}  (x^{-1}, x,y^{-1}x, x^{-1}y) \,\, \mbox{ and }  \,\,(x^{2^{n-2}+1}, x,y^{-1}x, x^{2^{n-2}+1}y)\end{equation}respectively. Now,  the automorphism of $G$ given by $(x,y) \mapsto (x^{-1}, x^{-1}y)$ in the first case, and the automorphism of $G$ given by $(x,y)\mapsto (x,xy^{-1})$ together with the braid transformation $\Phi_1$ in the second case show that the restrictions \eqref{resf} agree with \eqref{luz2} with $p=2^{n-2}$ and with $p=2$ respectively. Thus,  the action of $G$ on the surfaces lying in the strata defined by $\theta_{2^{n-2}}$ and $\theta_2$ extend to an action of $G_1$ and $G_2$ respectively, with signature $(0; 2,2,4,2^{n-1}).$

\s

{\it Claim.} The action of $G$ on $S$  does not extend to a group different from $G_1$ and $G_2.$

\s

Let us proceed by contradiction. Assume that the action of $G$ on some $S$ extends to an action of a group $H.$ By \cite[Theorem 1]{singerman2}, the group $H$ must have order $2^{n+1}$ and necessarily acts on $S$ with signature $(0; 2,2,4,2^{n-1}).$ Since $y^2$ is the unique involution of $G$, the group $H$ must contain an involution $z$ that does not belong to $G.$  Thus $$H \cong G \rtimes C_2 \,\, \mbox{ where } C_2 = \langle z \rangle$$and therefore $zxz=x^m$ where $m \in \{\pm 1, 2^{n-2} \pm1\}.$

\begin{enumerate}
\item[(a)] Assume $m=1.$ If we write $zyz=x^ky$ then the fact that $z^2=1$ shows that $k=0$ or $k=2^{n-2}.$ Observe that the former case is $H \cong G \times C_2$, but this group cannot be generated by two involutions and an element of order four. It follows that $zyz=y^{-1}$ and therefore $H=G_1.$
\s
\item[(b)] Assume $m=-1$ and write $zyz=x^ky.$ Observe that if $\hat{y}:=xy^{-1}$ then $z\hat{y}z=x^{k-2}\hat{y}$ and therefore we can assume that $k=0$ or $k=1.$ 

\s
\begin{enumerate}
\item The case $k=1$ must be disregarded since the group cannot be generated by two involutions and an element of order 4 with product of order $2^{n-1}.$ 

\s

\item If $k=0$ then $zxz=x^{-1}$ and $[z,y]=1.$ Note that $$\hat{z}:=x^{2^{n-3}}yz \,\, \mbox{ satisfies } \,\, \hat{z}^2=1, \, \hat{z}x\hat{z}=x \, \mbox{ and } \, \hat{z}y\hat{z}=y^{-1}.$$ Consequently,  $H$ is isomorphic to $G_1.$
\end{enumerate}
\s
\item[(c)] Assume $m=2^{n-2}-1$ and write $zyz=x^ky.$ The fact that $z^2=1$ implies that $k$ is even. Observe that  $$\hat{y}:=x^{2^{n-3}-1}y \, \mbox{ satisfies  } \, \hat{y}^4=1\, \mbox{ and } \,z\hat{y}z=x^{k+2}\hat{y}.$$ It follows that we can assume $k=2^{n-2}$ and  $zyz=y^{-1}.$ Thus, $H$ is isomorphic to  $G_2.$

\s

\item[(d)]  Assume $m=2^{n-2}+1$ and write $zyz=x^ky.$ The fact that $z^2=1$ implies that $4k \equiv 0 \mbox{ mod } 2^{n-1}$ and therefore $k \in \{0, 2^{n-2}, \pm 2^{n-3}\}$. The cases $k=\pm 2^{n-3}$ must be disregarded since they do not yield actions of order two.

\begin{enumerate}
\s

\item If $k=2^{n-2}$ then $zyz=y^{-1}.$ If $\hat{z}:=yz$ then $\hat{z}^2=1,$  $\hat{z}x\hat{z}=x^{2^{n-2}-1}$ and $\hat{z}y\hat{z}=y^{-1}.$ It then follows that $H$ is isomorphic to $G_2$.

\s

\item If $k=0$ then we let $\hat{y}:=xy$ and observe that $z\hat{y}z=(\hat{y})^{-1}$. It follows that this case agrees with case (a).
\end{enumerate}

 \end{enumerate}
The proof of the claim is done. The proof of the statement follows from the previous claim together with the fact that the signature $(0; 2,2,4,2^{n-1})$ is maximal.

\s

Finally, the  statement for each $S \in \mathscr{C}_{n,k}$ follows directly from the maximality of the signature $(2^{n-1}, 2^{n-k},4,4)$ for each $k \in \{2, \ldots, n-1\}.$

%
%

\section{Proof of Theorem \ref{des-fam-0}}

 Let $S$ belong to the family $\mathcal{F}_{n,0}.$ Following the proof of Theorem \ref{t-fam-sup}, the action of $G$ on each $S$ is represented by the ske $(x,xy, x^2)$ and the associated regular covering map $S \to S_G$ has exactly four branch points, all of them with $G$-isotropy group $\langle x^2 \rangle$. Observe that the intermediate regular covering map given by the action of $Z$ ramifies over exactly four branch values marked with $2.$ It follows from the Riemann-Hurwitz that the genus of $S_{Z}$ equals $2^{n-2}-1$ and therefore  \begin{equation}\label{r1} \dim \Pm(S \to S_{Z})=2^{n-2}.\end{equation} Similarly, for each $2 \leqslant j \leqslant n-2,$ the regular covering map given by the action of $H_j$ ramifies over two values marked with $2^{j-1}.$ Hence, the genus of $S_{H_j}$ equals $2^{n-j-1}$ and  \begin{equation}\label{r2} \dim \Pm(S_{H_j} \to S_{H_{j+1}})=2^{n-j-2}\end{equation} Now, the equalities \eqref{r1} and \eqref{r2}  together with the isogeny of Corollary \ref{isog-sup-gen} imply  $$2^{n-1}-1=1+\Sigma_{i=1}^3\dim \Pm(S_{N_i}\to S_G)+2^{n-2}+2(\Sigma_{j=2}^{n-2}2^{n-j-2})$$and, consequently, $S_{N_i}$ an elliptic curve for each $i$ and the result follows.

\s

Hereafter in this section we  denote by $z_1, z_2, z_3, z_4$ the ordered branch values of the regular covering map $\pi : S \to S_G.$

\s

 Let $S$ belong to the family $\mathcal{F}_{n,1}.$ Following the proof of Theorem \ref{t-fam-sup}, the action of $G$ on $S$ is represented by the ske $(xy,y,y,xy)$. We note that there are exactly $2^{n-2}$ branch points of $\pi$ over each $z_l$ and the $G$-isotropy group of them are $\langle x^iy\rangle$ where $i \in \mathbb{Z}_{2^{n-1}}$ is odd for $l=1,4$ and even for $l=2,3.$ The  fact that $Z \subset \langle x^iy \rangle$ for all $i$ allows us to see that the regular covering map given by the action of $Z$ on $S$ ramifies over $2^n$ values, each of them marked with 2. Thus, the Riemann-Hurwitz formula implies that the genus of $S_Z$ equals one and therefore $$\dim \Pm(S \to S_Z)=2^{n-1}.$$

In a similar way, observe that  $\langle x^iy \rangle \cap H_2$ equals $H_2=\langle y \rangle$ if $i=0$ and equals $\langle y^2 \rangle$ otherwise. It then follows that the regular covering map given by the action of $H_2$ on $S$ ramifies over four values marked with 4 and over $2^{n-1}-2$ values marked with 2. Hence, the Riemann-Hurwitz formula implies $S_{H_2}$ has genus zero. This fact implies, in turn, that the genera of $S_{N_2}$ (and then, by symmetry, the one of $JS_{N_3}$) and $S_{H_{j}}$ for each $j\geqslant 3$ are zero as well. The isogeny decomposition given in Corollary \ref{isog-sup-gen}  implies that   $JS_{N_1}$ is an elliptic curve, and the result follows

\s

 Let $S$ belong to the family $\mathcal{F}_{n,2}.$ Following the proof of Theorem \ref{t-fam-sup}, the action of $G$ on $S$ is represented by the ske $$(x, x^{p-1+2^{n-2} }, y,x^py) \,\, \mbox{ for some } p \in \mathbb{Z}_{2^{n-1}} \mbox{ that is even}.$$We observe that, independently of the stratum to which $S$ belongs (or independently of the choice of $p$) there are exactly two branch points of $\pi$ over each $z_l$ for $l=1, 2$ and the $G$-isotropy group of them equal $\langle x \rangle$, whereas there are exactly $2^{n-2}$ branch points of $\pi$ over each $z_l$ for $l=3,4$ and the $G$-isotropy group of have the form $\langle x^iy\rangle$ where $i$ is even.  

\s

 For each $j \in \{2, \ldots, n-1\}$ fixed, we consider the regular covering map $\pi_j$ given by the action of $H_j$ on $S$. Observe that the four branch values of $\pi$ with $G$-isotropy group $\langle x \rangle$ yield two branch values $\pi_j$ marked of  with $2^{j-1}$. In addition, the fact that
 $$\langle x^iy \rangle \cap H_j =  \langle x^iy \rangle \cong C_4 \mbox{ if } i \mbox{ is a multiple of } 2^{n-j}$$and $\langle x^iy \rangle \cap H_j=\langle y^2 \rangle$ otherwise, implies that among the $2^{n-1}$ branch values with $G$-isotropy group of the form $\langle x^i y\rangle$, there are exactly $2^{j}$ of them giving rise to four branch values of $\pi_j$ marked with 4, and the $2^{n-1}-2^{j}$ remaining ones give rise to $2^{n-j}-2$ branch values of $\pi_j$ marked with 2. In other words, the branch data of the action of $H_j$ on $S$ is  $$r(H_j)=(2^{j-1}, 2^{j-1}, 4 \stackrel{4}{\ldots} 4, 2,\stackrel{d}{\ldots},2) \, \mbox{ where } d=2^{n-j}-2.$$It then follows that the genus of $S_{H_j}$ is $2^{n-j-1}-1$ and, consequently,  $$\dim \Pm(S_{H_j} \to S_{H_{j+1}})=2^{n-j-2}.$$ By arguing similarly,  one obtains that$$r(Z)=(2, \stackrel{8}{\ldots}, 2), \,\, r({N_1})=(2^{n-1}, \stackrel{4}{\ldots}, 2^{n-1}, 2,2),  \,\, r(N_3)=(2^{n-2},2^{n-2},2,2)$$and therefore the genus of $S_{N_1}$ is zero, $JS_{N_3}$ is an elliptic curve and  $$\dim \Pm(S \to S_Z)=2^{n-1}.$$The isogeny decomposition of Corollary \ref{isog-sup-gen}  implies $$3 \cdot 2^{n-1}-1=0+\dim (JS_{N_2})+1+2^{n-1}+2(\Sigma_{j=2}^{n-2}2^{n-j-2})$$and therefore  $S_{N_2}$ has genus zero, and the result follows. Note that $$\Pm(S_{H_{n-2}} \to S_{H_{n-1}})=JS_{H_{n-2}} \mbox{ is an elliptic curve.}$$

\s

 Let $k \in \{2, \ldots, n-2\}$ and let $S$ belong to the family $\mathscr{C}_{n,k}.$ Following the proof of Theorem \ref{t-fam-sup}, the action of $G$ on $S$ is represented by a ske of the form $(g_1, g_2, y,xy)$ where $g_1$ generates $\langle x \rangle$ and $g_2$  generates the subgroup of $\langle x \rangle$ isomorphic to $C_{2^{n-k}}$. Without loss of generality, we can choose $g_2$ as $x^{2^{k-1}}$. 

The regular covering map given by the action of $Z$ on $S$ ramifies over $2^{n-1}+2^k+2$ values marked with 2. It follows that $$\dim JS_Z=2^{n-2}-2^{k-1} \,\, \mbox{ and therefore } \,\, \dim \Pm(S \to S_Z)=2^{n-1}.$$ In a similar way it can be checked that the genus of $S_{N_1}, S_{N_2}$ and $S_{N_3}$ equals zero.  

\s

In order to determine the genus of $S_{H_j}$ we have to proceed more carefully. Let $\pi_j$ denote the covering map given by the action of $H_j$ on $S$. Note that there are two branch 
 points of $\pi$ over $z_1$ and they give rise to one branch value of $\pi_j$ marked with $2^{j-1},$ and  there are $2^{n-2}$ branch points of $\pi$ over $z_4$ and they give rise to $2^{n-j-1}$ branch values of $\pi_j$  marked with $2$. Moreover, there are $2^{n-2}$ branch points of $\pi$ over $z_3$ and the $G$-isotropy group of them is $\langle x^py \rangle$ where $p \in \mathbb{Z}_{2^{n-1}}$ is even. We see that $2^{j-1}$ of these points produce two branch values of $\pi_j$ marked with 4 and the remaining $2^{n-2}-2^{j-1}$ ones produce $2^{n-j-1}-1$ branch values of $\pi_j$ marked with $2.$ Moreover,   there are $2^{k}$ branch 
 points of $\pi$ over $z_2$, the $G$-isotropy groups of them equal $\langle x^{2^{k-1}} \rangle$ and \begin{equation*}\label{ola}
\langle x^{2^{k-1}} \rangle \cap H_j = \left\{ \begin{array}{ll}
 \langle x^{2^{k-1}}\rangle \cong C_{2^{n-k}}  & \textrm{if $n-k \leqslant j-1$}\\
 \langle x^{2^{n-j}} \rangle \cong C_{2^{j-1}} & \textrm{if $n-k > j-1$}
  \end{array} \right.
\end{equation*}It follows that the $2^k$ branch points over $z_2$ yield $2^{n-j}$ branch values of $\pi_j$ marked with $2^{n-k}$ if $j \geqslant n-k+1,$ and yield $2^{k-1}$ branch values of $\pi_j$ marked with $2^{j-1}$ if $j < n-k+1$ 
\s

All the above says that the branch data of the action of $H_j$ on $S$ is 
 \begin{equation*}\label{ola}
r(H_j) = \left\{ \begin{array}{ll}
(2, \stackrel{d-1}{\ldots},2,4,4,2^{j-1}, 2^{n-k}, \stackrel{d}{\ldots}, 2^{n-k})  & \textrm{if $j \geqslant n-k+1$}\\
(2, \stackrel{d-1}{\ldots},2, 4,4,2^{j-1}, \stackrel{d'}{\ldots}, 2^{j-1}) & \textrm{if $j < n-k+1$}
  \end{array} \right.
\end{equation*}where $d=2^{n-j}$ and $d'=2^{k-1}+1$. It then follows from Riemann Hurwitz formula that $$\dim JS_{H_j}=2^{n-j-1}-2^{k-2} \, \, \mbox{ if } \,\, j < n-k+1$$and zero otherwise. The result follows from the isogeny given in  Corollary \ref{isog-sup-gen}.

\section{The classical quaternion group} \label{TN}

In this section we collect results  concerning abelian varieties $A$ and compact Riemann surfaces $S$ with a group of automorphisms $G$ isomorphic to the quaternion group of order eight $Q(8)$. For the sake of brevity we shall avoid some details, but they can be deduced by arguing as we did to prove the previous results.

\s

The isotypical decomposition of $A$ with respect to $G$ (which agrees with the group algebra decompositions) is given by
\begin{equation}\label{inc}
A \sim A_{G} \times \Pi_{i=1}^3   \Pm(A_{N_i} /A_{G}) \times \Pm(A/ A_{Z}).\end{equation} 

If the analytic representation of the action of $G$ on $A$ decomposes as  $$\rho_a \simeq \oplus_{j=1}^4 a_j \chi_j \oplus b\Theta_1$$where $a_i, b \geqslant 0$ are integers, then the dimension of the factors arising in \eqref{inc} are $$\dim A_G=a_1, \,\, \dim \Pm(A_{N_i}/A_{G})=a_{i+1} \,\, \mbox{ and }\,\, \dim \Pm(A/A_{Z})=2b$$

Moreover, the following statements are equivalent:

\begin{enumerate}
\item The isotypical decomposition of $A$ with respect to $G$ is trivial.
\item The dimension of $A_Z$ is zero.
\item The dimension of $A_{K}$ is zero, for each nontrivial subgroup $K$ of $G.$
\item The integers $a_i$ equal zero, for all $i.$ 
\item 1 is not an eigenvalue of $\rho_a(g)$ for each nontrivial $g \in G.$
\end{enumerate}

The previous results applied to $A=JS$ provide analogous results for Jacobians, as in the general case.
\s

 There exist precisely three complex one-dimensional families of Riemann surfaces endowed with a group of automorphisms $G$ isomorphic to  $Q(8)$. 

\subsection*{(1)} The family $\mathcal{F}_{3,0}$ consisting of those Riemann surfaces $S$ of genus three with $G$ acting on them with signature $(1;2).$ This family consists of only one equisymmetric stratum and the action of $G$ on each $S$ extends to an action of \begin{equation}\label{obligado}\langle a,b,c : a^2,b^2, c^4, bcbc^3, acac^3, abac^2b\rangle =  (C_4 \times C_2) \rtimes C_2\end{equation}with signature $(0; 2,2,2,4).$ Moreover, up to finitely many exceptional surfaces, the full automorphism group of $S$  is isomorphic to \eqref{obligado}. See also \cite[Table 5, 3.ad.2]{Brou}. With respect to this presentation, $G$ is isomorphic to $\langle x:=ca, y:=ba\rangle.$

The isotypical decomposition of $JS$ with respect to $G$ turns into the canonical decomposition $$JS \sim E \times P$$where $E \sim JS_G$ is an elliptic curve and $P \sim \mbox{Prym}(S \to S_G)$ an abelian surface. However,  $JS$ decomposes further. Indeed, by considering group algebra decompositions with respect to \eqref{obligado}, it decomposes as $$JS \sim E_1^2 \times E_2$$where $E_1\sim JS_{\langle a \rangle}$ and $E_2\sim JS_{\langle ab \rangle}$ are elliptic curves.

\subsection*{(2)}The family $\mathcal{F}_{3,1}=\mathcal{F}_{3,2}$ consisting of those   Riemann surfaces $S$ of genus five with $G$ acting on them with signature $(0; 4,4,4,4).$ This family consists of only one equisymmetric stratum and the action of $G$ on each $S$ extends to an action of the group isomorphic to $(\mathbb{D}_4 \times C_2) \rtimes C_2$ presented as \begin{equation}\label{caldo}\langle r,s,a,b : r^4,s^2, a^2, b^2, (sr)^2, arar^{-1}, (as)^2, 
brbr^{-1}, bsb(sra)^{-1}, bab(ar^2)^{-1} \rangle \end{equation}with signature $(0; 2,2,2,4)$. 
Moreover, up to finitely many exceptional surfaces, the full automorphism group of $S$  is isomorphic to \eqref{caldo}. With respect to the preceding presentation, $G$ is isomorphic to $\langle x:=ra, y:=rb\rangle.$  See also \cite[Lemma 8]{CI45}.

The isotypical decomposition of $JS$ with respect to $G$ is $$JS \sim E \times P'$$where $E \sim JS_Z$ is an elliptic curve and $P' \sim \mbox{Prym}(S \to S_Z)$ an abelian fourfold. However,  by considering group algebra decompositions of $JS$ with respect to  \eqref{caldo}, the Jacobian decomposes further as  $$JS \sim E_1 \times E_2^4$$where $E_1\sim JS_{\langle r \rangle}$ and $E_2\sim JS_{\langle s,a \rangle}$ are elliptic curves.

\subsection*{(3)} The family $\mathscr{C}_{3,2}$ consisting of those Riemann surfaces $S$ of genus four with $G$ acting on them with signature $(0; 2,4,4,4).$ This family consists of only one equisymmetric stratum and, up to finitely many exceptional surfaces, the full automorphisms group of $S$ is isomorphic to $G$. See also \cite[Lemma 3.1]{CIg4}.  These surfaces are hyperelliptic and their isotypical  decompositions are trivial.

\section{Proof of Proposition \ref{malg}}
Let $n \geqslant 3$ be an integer and let $S$ be a  Riemann surface lying in the family $\mathscr{C}_{n,n-1}.$ The two-fold regular covering map $\pi_Z:S  \to \mathbb{P}^1$ given by the action of $Z=\langle y^2 \rangle$ on $S$ ramifies over $2^n+2$ branch values. If we denote these values by \begin{equation}\label{val}\delta_1, \delta_2, \alpha_1, \ldots, \alpha_{2^{n-2}}, \alpha'_1, \ldots, \alpha'_{2^{n-2}}, \beta_1, \ldots, \beta_{2^{n-2}}, \gamma_{1}, \ldots, \gamma_{2^{n-2}}\end{equation}then $S$ is isomorphic to the normalization of the singular affine algebraic curve \begin{equation} \label{finall} Y^2=(X-\delta_1)(X-\delta_2)\Pi_{j=1}^{2^{n-2}}(X-\alpha_j)(X-\alpha_j')(X-\beta_j)(X-\gamma_j).\end{equation}

Note that $$K=G/Z \cong \mathbb{D}_{2^{n-2}}=\langle r, s : r^{2^{n-2}}, s^2, (sr)^2 \rangle $$acts on $S_Z \cong \mathbb{P}^1$ in such a way that the corresponding quotient is isomorphic to $S_G.$ Without loss of generality, we can assume $$r(z)=\omega z \,\, \mbox{ where } \omega=\mbox{exp}(\tfrac{2 \pi i}{2^{n-2}}).$$The fact that the branch values \eqref{val} are setwise invariant under the action of $r$ implies that we can assume that there exist nonzero complex numbers $\lambda_2, \lambda_3, \lambda_4$ different from 1 and satisfying that $\lambda_i \notin \{\lambda_l\omega^{k}: 1 \leqslant k \leqslant 2^{n-2}\}$ for each $i \neq l$, such that  \begin{equation}\label{az}\delta_1=0, \delta_2=\infty, \alpha_j=\omega^j, \alpha'_j=\lambda_2\omega^j, \beta_j=\lambda_3\omega^j \,\, \mbox{ and }\,\, \gamma_j=\lambda_4\omega^j\end{equation}for each $j \in \{1, \ldots, 2^{n-2}\}.$ The covering map induced by $\langle r \rangle$ can be assumed to be $z \mapsto z^{2^{n-2}}$ and therefore the composed covering map   $S \to S_Z \to (S_Z)_{\langle r \rangle}$ ramifies over $\infty$ and $0$ marked with $2^{n-1},$ and over $$1,\,\, \lambda_2^{2^{n-2}},\lambda_3^{2^{n-2}}  \mbox{ and }\lambda_4^{2^{n-2}}$$marked with 2.  Note that $s$ induces an involution $\tilde{s}$ of $(S_Z)_{\langle r \rangle}$ in such a way that the corresponding quotient is isomorphic to $S_G.$ Thus, $\tilde{s}$ acts with two orbits of length two and with two fixed points. Without loss of generality, we can assume $$\tilde{s}(z)=\lambda_2^{2^{n-2}}\tfrac{1}{z}$$and therefore the orbits are $\{\infty, 0\}$, $\{1, \lambda_2^{2^{n-2}}\}$ and the fixed points are $\lambda_3^{2^{n-2}}$ and $\lambda_4^{2^{n-2}}.$ It follows that $$ \lambda_2^{2^{n-2}}=\lambda_3^{2^{n-1}}=\lambda_4^{2^{n-1}} \mbox{ showing that }\lambda_4^{2^{n-2}}=- \lambda_3^{2^{n-2}}.$$

If we let $t:=\lambda_2^{2^{n-2}}$ then $\lambda_3^{2^{n-2}}=\sqrt{t}$ and $\lambda_4^{2^{n-2}}=-\sqrt{t}.$ Note that $t \in \mathbb{C}-\{0,1\}.$ Now, we replace \eqref{az} in  \eqref{finall} to obtain  that $S$ is isomorphic to the normalization of \begin{equation}\label{neg}Y^2=X(X^{2^{n-2}}-1)(X^{2^{n-2}}-t)(X^{2^{n-1}}-t)\end{equation}as desired, where the factor $X-\delta_2$ has been disregarded since $\delta_2=\infty.$ Finally, if  $$\mathbf{x}(X,Y) = (\xi^2 X, \xi Y) \,\, \mbox{ and } \,\, \mathbf{y}(X,Y) = (\tfrac{\lambda}{X}, \eta t \tfrac{Y}{X^{2^{n-1}+1}})$$where $\xi=\exppp(\tfrac{2 \pi i}{2^{n-1}}),$ $\eta^2=-\lambda$ and $\lambda^{2^{n-2}}=t$ then it is a routine computation to verify that $\mathbf{x}$ and $\mathbf{y}$ are automorphisms of \eqref{neg} and generate a group isomorphic to $Q(2^n).$

\section{Proof of Theorems \ref{ppavf1}, \ref{ppavf2} and  \ref{ppavf3} and Proposition \ref{especial}}

Let $S$ be a compact Riemann surface of genus $g \geqslant 2.$ Consider the Jacobian variety $JS \cong \mathbb{C}^g/\Lambda$ and its full (polarization-preserving) automorphism group  $\mbox{Aut}(JS).$ Every automorphism of  $S$ induces a unique automorphism of $JS$. In fact $$[\mbox{Aut}(JS): \mbox{Aut}(S)] \in \{1,2\}$$according to whether or not $S$ is hyperelliptic; moreover, in the latter case $$ \mbox{Aut}(JS)/ \mbox{Aut}(S)=\{\pm  \mbox{id}\}.$$

Once a symplectic basis of $\Lambda=H_{1}(S, \mathbb{Z})$ is fixed, there is an isomorphism  $$\mbox{Aut}(JS) \cong \Sigma_S:=\{ R \in \Sp(2g, \mathbb{Z}) : R \cdot Z_S=Z_S\},$$where $(I_g \, Z_S)$ is the period matrix of $JS$. Two different basis induce different but equivalent matrices $Z_S$ and conjugate subgroups $\Sigma_S.$ Thus, we obtain a well-defined analytic submanifold of $\mathscr{H}_g$ given by $$\mathscr{S}_S:=\{Z\in \mathscr{H}_g: R \cdot Z=Z \text{ for all } R \in \Sigma_S\}$$ whose points are matrices representing principally polarized abelian varieties admitting an action  equivalent to the one of $\mbox{Aut}(S)$ in the symplectic group. Observe that $-I_{g} \in \mathscr{S}_S.$ 

\s

Following the results of  \S\ref{TN}, the family $\mathcal{F}_{3,0}$ consists of one  stratum and, up to finitely many surfaces, the full automorphism group of its members is isomorphic to $$\mathscr{G}=\langle a,b,c : a^2,b^2, c^4, bcbc^3, acac^3, abac^2b\rangle =  (C_4 \times C_2) \rtimes C_2$$and acts on then with signature $(0; 2,2,2,4).$ The action is represented by the ske $$\theta=(a,b,abc^{-1},c^{-1}).$$

We apply the results on adapted hyperbolic polygons, geometric generators and symplectic representations obtained in \cite{poly} as well as the algorithms  programed  in the same article to represent symplectically the group $\mathscr{G}$ starting from the action $\theta.$ Explicitly, the correspondence$$a \mapsto \left(\begin{smallmatrix}
0&0&0&0&-1&0\\
-1&0&0&1&0&0\\
0&0&-1&0&0&0\\
0&1&0&0&-1&0\\
-1&0&0&0&0&0\\
0&0&0&0&0&-1
\end{smallmatrix}\right) \,\,\,\,\,\,\,\, b \mapsto \left(\begin{smallmatrix}
-1&1&0&0&-1&1\\
-1&0&-1&1&0&1\\
1&0&0&-1&-1&0\\
0&1&0&-1&-1&1\\
-1&0&-1&1&0&0\\
0&1 &0&0&-1&0
\end{smallmatrix}\right) \,\,\,\,\,\,\,\, c \mapsto \left(\begin{smallmatrix}
0 & 0 & 0 & -1 &-1&1\\
0&0&0&0&-1&1\\
0&0&0&0&0&-1\\
1&0&0&-1&-1&0\\
-1&1&0&0&-1&1\\
0&1&1&-1&-1&0 
\end{smallmatrix}\right)$$defines an isomorphism between $\mathscr{G}$ and a subgroup $\Sigma_{3,0}$ of $\mbox{Sp}(6, \mathbb{Z}).$

\s

Now, if $Z \in \mbox{GL}(3, \mathbb{C})$ is symmetric and satisfy$$R \cdot Z = (AZ+B)(CZ+D)^{-1}=Z \,\, \mbox{ for each }\,\, R=  \left( \begin{smallmatrix}
A & B \\
C & D
\end{smallmatrix} \right) \in \Sigma_{3,0}$$then, after some computations, we obtain that $Z$ equals $$Z_t=\left( \begin{smallmatrix}
1+i+\tfrac{1-i}{2}t & 1-it & \tfrac{1+i}{2}t \\
1-it & 1-(1+i)t & t \\
\tfrac{1+i}{2}t & t & 1+\tfrac{i-1}{2}t
\end{smallmatrix} \right)$$for some $t \in \mathbb{C},$ where $i=\sqrt{-1}.$ If $\pi_3: \mathscr{H}_3 \to \mathcal{A}_3$ is the canonical projection (see \S\ref{esp}) then  $$\mathscr{Y}_3:=\pi_3(\{Z_t : t \in \mathbb{C}\}\cap \mathscr{H}_3) \subset \mbox{Sing}(A_3)$$is a complex one-dimensional family of principally polarized abelian threefolds  whose members, by construction,  admit an action of $\mathscr{G}$ with rational representation $\mathscr{G} \to \Sigma_{3,0}.$  Clearly, this family contains the image under the Torelli map of the family $\mathcal{F}_{3,0}$

\s

The fact that all the members of $\mathscr{Y}_3$ share the same rational representation of $\mathscr{G}$ implies that the analytic representation $\rho_a$ of $G=\langle ca,ba\rangle \cong Q(8)$ on each member of $\mathscr{Y}_3$ is equivalent to the one of  the members of $\mathcal{F}_{3,0}$.  Thus, $\rho_a \simeq \chi_1 \oplus \Theta_1$. 

\s

The arguments employed to prove Proposition \ref{pp1} together with the facts discussed in \S\ref{TN} show that the members of $\mathscr{Y}_3$ decompose isogenously as the product of an elliptic curve and an abelian surface. This completes the proof of Theorem \ref{ppavf1}.

\s

To prove Theorem \ref{ppavf2} we proceed analogously, but taking into account that the  the full automorphism group of a generic member of the family $\mathcal{F}_{3,1}=\mathcal{F}_{3,2}$ 
 is isomorphic to$$\mathscr{G}'=\langle r,s,a,b : r^4,s^2, a^2, b^2, (sr)^2, arar^{-1}, (as)^2, 
brbr^{-1}, bsb(sra)^{-1}, bab(ar^2)^{-1} \rangle$$and that the image $\Sigma_{3,1} \leqslant \mbox{Sp}(10, \mathbb{Z})$ of its symplectic representation  is generated by $$\left(\begin{smallmatrix}
-2& -2&  1&  1&  1& -3&  0&  1&  5&  4\\
 1&  1&  0 &-1 &-1 & 4 & 0 & 0 &-1& -3\\
-1& -1 & 1&  0& -1 & 1 & 0&  0 & 0 & 2\\
 0 & 1 & 0 & 0  &0 & 1 &-1 & 0 & 0 & 0\\
 0 & 0 & 1&  0 &-1 & 2 & 0 & 2 &-3 & 2\\
0 & 0 & 0 & 0 & 0 & 0 & 0 & 1 &-1 & 1\\
0 & 0  &0 & 0 & 0 & 0 & 0 & 1 &-2 & 1\\
0 & 0 & 0 & 0 & 0 &-1 & 0 &-1 & 0& -1\\
0 & 0 & 0 & 0 & 0 &-1 & 1 & 0 & 0 & 0\\
 0 & 0 & 0&  0 & 0 & 1&  0 & 1 & 1&  0
\end{smallmatrix}\right) \,\,\,\, \left(\begin{smallmatrix}
-2& -2 & 1&  1 & 1 & 0 &-3 & 1 & 5  &4\\
 2  &2& -1 &-1 & 0 & 3 & 0 & 0 &-1 &-5\\
 1 & 1& -1 & 0 & 1& -1 & 0 & 0&  0 &-2\\
-1& -2 & 1 & 0 & 1 &-5 & 1 & 0 & 0 & 3\\
 1 & 1 & 0 & 0 & 0 &-4 & 5  &2 &-3 & 0\\
 0  &0 & 0 & 0 & 0 &-2 & 2 & 1& -1 & 1\\
 0 & 0 & 0 & 0 & 0 &-2 & 2 & 1& -2 & 1\\
 0 & 0 & 0 & 0 & 0 & 1 &-1 &-1 & 1 & 0\\
 0 & 0&  0 & 0 & 0 & 1& -1 & 0 & 0 & 0\\
0 & 0&  0 & 0 & 0 & 1&  0 & 1  &1 & 0
\end{smallmatrix}\right)\,\,\,\, \left(\begin{smallmatrix}
-1 &-2&  0  &0 & 2 & 0 &-2 & 0 & 6 & 4\\
 0 & 0 & 0 & 0 &-1&  2 & 0 & 0 & 0 &-1\\
 0 & 0 &-1 & 0 & 0 & 0 & 0 & 0 & 1 & 0\\
-2 &-3  &1 & 1 & 1 &-6 & 0& -1&  0&  5\\
 0 &-1&  0&  0 & 0 &-4 & 1 & 0& -5 & 0\\
 0 & 0 & 0 & 0 & 0 &-1 & 0 & 0 &-2  &0\\
 0 & 0 & 0 & 0 & 0 &-2 & 0 & 0 &-3 &-1\\
 0 & 0 & 0 & 0  &0 & 0 & 0 &-1 & 1 & 0\\
 0 & 0 & 0 & 0&  0 & 0 & 0 & 0  &1 & 0\\
 0 & 0 & 0 & 0&  0 & 2 &-1 & 0  &1 & 0
\end{smallmatrix}\right)$$

\s

To guarantee the existence of the family $\mathscr{Y}_4$ and to prove the statements of Theorem \ref{ppavf3} we also proceed as in the preceding cases. Since we are enable to solve the corresponding system of equations $R \cdot Z=Z$ (it requires solving a system of several fairly complicated quadratic  equations in ten variables, even more complicated than in the apparently more challenging case of $\mathscr{Y}_5$), we are not in position to provide an explicit description of the corresponding period matrices and, consequently, we cannot deduce from that its dimension. However, the complex dimension of the family $\mathscr{Y}_4$ equals 1, as it was computed in \cite{frediani} (the family $\mathscr{Y}_4$ is labeled there as the case (36) in Table 1).

\s

Finally, we apply again the results of \cite{poly} to determine a symplectic representation  of \begin{equation*}\label{pan}\mbox{Aut}(X_4) \cong Q\mathbb{D}_{16} =\langle u,v : u^{16}, v^2, vuvu^{-7}\rangle \end{equation*}according to its action on $X_4.$ The image $\Sigma_0$ of this representation $Q\mathbb{D}_{16} \to \Sigma_0$ is generated by the matrices$$A=\left(\begin{smallmatrix}
 0 & 0 & 0 & 0 &-1 & 0 & 0 & 0 \\
 0 & 1 & 1 & 1 & 1 &-2 &-1 & 0\\
 0 & 0 & 1 & 1 & 0 & 1 &-2&  0\\
0 & 0  &0 & 1  &0 & 0  &1 &-1\\
1&  1  &1  &1 &-1 &-1 &-1 & 0\\
0 & 1 & 1 & 1 & 0& -1 &-1 & 0\\
0 & 0 & 1  &1&  0 & 0 &-1 & 0\\
 0 & 0  &0 & 1 & 0 & 0 & 0 & 0
 \end{smallmatrix}\right)
        \,\,\,\,\,\,\, B=\left(\begin{smallmatrix}
 0& -1 &-1 &-1 &-1  &1 & 1 & 1\\
 0 & 0&  0 & 0 & 1 &-1 & 0  &1\\
-1& -1 &-1 & 0 & 1 & 0 & 2& -1\\
-1& -1 & 0 & 0  &1 & 1 &-1 & 0\\
0 &-1 &-1& -1 & 0 & 0 & 1  &1\\
-1& -1 &-1 &-1 & 1 & 0 & 1&  1\\
-1 &-1& -1 & 0 & 1 & 0 & 1 & 0\\
-1 &-1 & 0  &0 & 1  &0 & 0 & 0
\end{smallmatrix}\right)$$Now, after some computations, one sees that the equations$$A \cdot Z= B \cdot Z = Z=Z^t, \,Z \in \mbox{GL}(4, \mathbb{C})$$admit a unique solution $Z_0$ in $\mathscr{H}_4.$ In other words, there exists a unique principally polarized abelian fourfold $A$ endowed with an action of $Q\mathbb{D}_{16}$ with rational representation equivalent to $Q\mathbb{D}_{16} \to \Sigma_0$. Clearly, $A=JX_4$ and $Z_0$ is its period matrix; this matrix is the one given in the statement of Proposition \ref{especial} and the proof is done.


\begin{thebibliography}{9}
\bibitem{CI45} 
{\sc G. Bartolini, A. F. Costa, and M. Izquierdo,} 
{\em On the orbifold structure of the moduli space of Riemann surfaces of genera four and five.}
Rev. R. Acad. Cienc. Exactas Fis. Nat. Ser. A Mat. RACSAM {\bf 108} (2014), no. 2, 769--793. 

\bibitem{poly}
{\sc A. Behn, R. E. Rodr\'iguez and A. M. Rojas,} 
{\em Adapted hyperbolic polygons and symplectic representations for group actions on Riemann surfaces,} 
J. Pure Appl. Algebra {\bf 217} (2013), no. 3, 409--426.

\bibitem{Brou}
{\sc{S. A. Broughton,}}  {\em{Classifying finite groups actions on surfaces of low genus,}} J. Pure Appl. Algebra {\bf 69} (1990), no. 3, 233--270. 

\bibitem{b}
{\sc S. A. Broughton},
 { \em The equisymmetric stratification of the moduli space and the Krull dimension of mapping class groups,} 
 Topology Appl. {\bf 37} (1990), no. 2, 101--113.
 
 \bibitem{Betal}
{\sc E. Bujalance, J. J. Etayo and E. Martinez,} {\em Automorphisms groups of hyperelliptic Riemann surfaces,} Kodai Math. J. {\bf 10} (1987)
174--181.
 
 
\bibitem{d1}
{\sc A. Carocca, S. Recillas and R. E. Rodr\'iguez}, {\em Dihedral groups acting on Jacobians,} Contemp. Math. {\bf 311} (2011), 41--77.

\bibitem{CR-H}
{\sc A. Carocca and R. E. Rodr\'iguez,}
{\em Hecke algebras acting on abelian varieties},
J. Pure Appl. Algebra {\bf 222} (2018), no. 9, 2626--2647.



\bibitem{CR}
{\sc A. Carocca and R. E. Rodr\'iguez,}
{\em Jacobians with group actions and rational idempotents.}
J. Algebra \textbf{306} (2006), no. 2, 322--343.

\bibitem{CIg4}
{\sc A. F. Costa and M. Izquierdo,}
{\em Equisymmetric strata of the singular locus of the moduli space of Riemann surfaces of genus 4,} Geometry of Riemann surfaces, 120--138, London Math. Soc. Lecture Note Ser., 368, Cambridge Univ. Press, Cambridge, 2010. 


\bibitem{frediani} 
{\sc P. Frediani, A. Ghigi and M. Penegini}, 
{\em Shimura varieties in the {T}orelli locus via {G}alois coverings}, 
 Int. Math. Res. Not. {\bf 20}  (2015) 10595--10623.
 
 \bibitem{Grot}
{\sc A. Grothendieck,} {\em Esquisse d'un Programme}. In Geometric Galois Actions 1, Around Grothendieck's Esquisse d'un Programme. ed. Schneps, L. and Lochak, P., London Math. Soc. Lecture Note Ser. {\bf 242} (1997), 5--48.

 
\bibitem{Harvey}
{\sc J. Harvey,}
{\em On branch loci in Teichm\"{u}ller space},
Trans. Amer. Math. Soc. {\bf 153} (1971), 387--399.

\bibitem{nos}
{\sc R. A. Hidalgo, L. Jim\'enez, S. Quispe and S. Reyes-Carocca,} {\em Quasiplatonic curves with symmetry group $\mathbb{Z}_2^2 \rtimes \mathbb{Z}_m$ are definable over $\mathbb{Q}$,} Bull. London Math. Soc. {\bf 49} (2017) 165--183.


\bibitem{HQ}
{\sc{R. A. Hidalgo and S. Quispe,}}  {\em{Regular dessins d'enfants with dicyclic group of automorphisms,}} J. of Pure Appl. Algebra {\bf 224} no. 5 (2020).


\bibitem{KS}
{\sc{S. Kallel and D. Sjerve,}}  {\em{Genus zero actions on Riemann surfaces,}} Kyushu J. Math. {\bf 55} (2001), no. 1, 141--164.

\bibitem{LR} 
{\sc H. Lange and S. Recillas,}
{\em Abelian varieties with group action}.
J. reine angew. Math. {\bf 575} (2004), 135--155.
 
\bibitem{oort}
{\sc F. Oort},
 { \em Singularities of coarse moduli schemes,} S\'em. Dubriel {\bf 16} (1976).
 

\bibitem{PA}
{\sc J. Paulhus and A. M. Rojas,}
{ \em Completely decomposable Jacobian varieties in new genera},  Experimental Mathematics {\bf 26} (2017), no. 4, 430--445.

 \bibitem{jpaa}
{\sc S. Reyes-Carocca,} {\em On the one-dimensional family of Riemann surfaces of genus $q$ with $4q$ automorphisms}, J. Pure Appl. Algebra 223, no. {\bf 5} (2019), 2123--2144.



\bibitem{pisa}
{\sc  S. Reyes-Carocca and R. E. Rodr\'iguez},  
{\em A generalisation of Kani-Rosen decomposition theorem for Jacobian varieties,} Ann. Sc. Norm. Super. Pisa Cl. Sci. (5) {\bf 19} (2019), no. 2, 705--722. 


\bibitem{RCR}
{\sc S. Reyes-Carocca and R. E. Rodr\'iguez,} 
{\em On Jacobians with group action and coverings}, 
Math. Z. (2020) {\bf 294}, 209--227.

\bibitem{d2}
{\sc J. Ries}, {\em The Prym variety for a cyclic unramified cover of a hyperelliptic curves}, J. Reine Angew. Math. \textbf{340} (1983) 59--69.

\bibitem{anita-ibero}
{\sc A. M.  Rojas,}
{\em Group actions on Jacobian varieties}, 
Rev. Mat. Iber. {\bf 23}  (2007),  no. 2, 397--420. 

\bibitem{d4}
{\sc A. S\'anchez-Arg\'aez}, {\em Actions of the group $A_5$ in Jacobian varieties}, Aportaciones Mat. Comun. {\bf 25}, Soc. Mat. Mexicana, M\'exico (1999), 99--108.


\bibitem{Serre} 
{\sc J-P. Serre,}
{\em  Linear Representations of finite groups,}
Graduate text in Mathematics {\bf 42}, 1996. 

\bibitem{singerman2}
{\sc D. Singerman}, 
{\em Finitely maximal Fuchsian groups}, 
J. London Math. Soc. (2)  {\bf 6}, (1972), 29--38.

\bibitem{singerman}
{\sc D. Singerman}, {\em Subgroups of Fuchsian groups and finite permutation groups}, Bull. London Math. Soc.  {\bf 2}, (1970), 319--323.

\bibitem{W}
{\sc J. Wolfart,} {\em Triangle groups and Jacobians of CM type.} Frankfurt a.M., 2000. www.math.uni-frankfurt.de/$\sim$wolfart/Artikel/jac.pdf

\end{thebibliography}
\end{document}